\documentclass{amsart}
\usepackage{amsfonts}
\usepackage{amsmath}
\usepackage{amscd}
\usepackage{amssymb}
\usepackage{latexsym}
\usepackage{amsthm}

\newcommand \aff{\operatorname{Sp}} 
\newcommand \after{\circ}

\newcommand \into{\hookrightarrow}

\newcommand \inverse[2]{{#1^{-1}(#2)}}
\newcommand \iso{\cong}
\newcommand \loc{{\mathcal {O}}}
\newcommand \map[1]{{\newcommand{\tmpprop}{#1q} \newcommand{\emptyprop}{q} \if\tmpprop\emptyprop \to\else \xrightarrow{\phantom{i}#1\phantom{i}}\fi}} 
\newcommand \maxim{\mathfrak m}
\newcommand \nat{\mathbb N}
\newcommand \norm[1]{\left|#1\right|}
\newcommand \onto{\twoheadrightarrow}
\newcommand \operate[2]{{\operatorname{#1}(#2)}}

\newcommand \pr{\mathfrak p}

\newcommand \range [2]{#1,\dots,#2}

\newcommand \restrict [2]{\left.#1\right|{}_{{#2}}}
\newcommand \rij[2]{(#1_1,\dots,#1_{#2})}
\newcommand \scp[2]{#1\langle#2\rangle}
\newcommand \set[2]{\left\{\,#1\mid #2\,\right\}} 

\newcommand \tensor{\otimes}

\newcommand \op\operatorname

\newcommand\En{\bigwedge}
\newcommand\en{\wedge}
\newcommand\niet{\neg}
\newcommand\Of{\bigvee}
\newcommand\of{\vee}

 \newcommand\commtrianglefront[7][]{%
\begin{equation}
{\newcommand{\tmpprop}{#1q} \newcommand{\emptyprop}{q}
\if\tmpprop\emptyprop \relax\else \label{#1}\fi}
\begin{aligned}%
\mbox{
\begin{picture}(120,80)%
\put(55,70){\vector(-1,-2){30}}
\put(65,70){\vector(1,-2){30}}
\put(30,5){\vector(1,0){60}}
\put(60,75){\makebox(0,0)[c]{$#2$}}
\put(25,5){\makebox(0,0)[r]{$#4$}}
\put(95,5){\makebox(0,0)[l]{$#6$}}
\put(60,0){\makebox(0,0)[c]{$#5$}}
\put(37,43){\makebox(0,0)[r]{$#3$}}
\put(83,43){\makebox(0,0)[l]{$#7$}}
\end{picture}}
\end{aligned}
\end{equation}}

\newcommand\commtriangleback[7][]{%
\begin{equation}
{\newcommand{\tmpprop}{#1q} \newcommand{\emptyprop}{q}
\if\tmpprop\emptyprop \relax\else \label{#1}\fi}
\begin{aligned}%
\mbox{
\begin{picture}(120,80)%
\put(55,70){\vector(-1,-2){30}}
\put(65,70){\vector(1,-2){30}}
\put(30,5){\vector(1,0){60}}
\put(60,75){\makebox(0,0)[c]{$#2$}}
\put(25,5){\makebox(0,0)[r]{$#6$}}
\put(95,5){\makebox(0,0)[l]{$#4$}}
\put(60,0){\makebox(0,0)[c]{$#5$}}
\put(37,43){\makebox(0,0)[r]{$#7$}}
\put(83,43){\makebox(0,0)[l]{$#3$}}
\end{picture}}
\end{aligned}
\end{equation}}

\newcommand\commdiagram[9][]{%
\begin{equation}
{\newcommand{\tmpprop}{#1q} \newcommand{\emptyprop}{q}
\if\tmpprop\emptyprop \relax\else \label{#1}\fi}
\begin{aligned}%
\mbox{
\begin{picture}(130,90)%
\put(120,70){\vector( 0,-1){50}}%
\put(10,80){\vector( 1, 0){100}}%
\put(0,70){\vector( 0,-1){50}}%
\put(10,10){\vector( 1, 0){100}}%
\put(115,80){\makebox(0,0)[l]{$#4$}}%
\put(5,80){\makebox(0,0)[r]{$#2$}}%
\put(115,10){\makebox(0,0)[l]{$#9$}}%
\put(5,10){\makebox(0,0)[r]{$#7$}}%
\put(-3,50){\makebox(0,0)[r]{$#5$}}
\put(123,50){\makebox(0,0)[l]{$#6$}}
\put(60,3){\makebox(0,0)[c]{$#8$}}
\put(60,88){\makebox(0,0)[c]{$#3$}}
\end{picture}}
\end{aligned}
\end{equation}}

\newcommand\name[1]{{\sc#1}}

\newcommand\refine{{\prec}}
\newcommand\regsize[1]{{\mathstrut\smash{#1}}}
\newcommand\Berk[1]{{\mathbb{#1}}}
\newcommand\anDlang{{\mathcal L_{\operatorname{an}}^{\mathbf D}}} 
\newcommand\anlang{{\mathcal L_{\operatorname{an}}}} 
\newcommand\adm[1]{{#1^\circ}} 
 \newcommand\tuple[1]{{\mathbf{#1}}}

\newcommand\ERS{Embedded Resolution of Singularities}
\newcommand\QE{Quantifier Elimination}
\newcommand\rav{rigid analytic variety}
\newcommand\ravs{rigid analytic varieties}

\newcommand\VE{Vo\^ute Etoil\'ee}
\newcommand\etoile{\'etoile}
\newcommand\Loj{\L ojasiewicz}

\newcounter{lst}
\newcounter{stp}

	{\begin{enumerate}\setcounter{enumi}{\value{lst}}}
	{\setcounter{lst}{\value{enumi}}\end{enumerate}}
\newenvironment{step}
	{\relax\par\smallskip\noindent\addtocounter{stp}1\emph{Step \thestp}.}
	{}
\newenvironment{interlude}
	{\relax\par\smallskip\noindent\emph{Interlude}.}
	{}

\swapnumbers
\theoremstyle{plain}
\newtheorem{theorem}[subsection]{Theorem}
\newtheorem{corollary}[subsection]{Corollary}
\newtheorem{lemma}[subsection]{Lemma}
\newtheorem{proposition}[subsection]{Proposition}

\theoremstyle{definition}
\newtheorem{definition}[subsection]{Definition}

\newtheorem{example}[subsection]{Example}

\theoremstyle{remark}

\newtheorem{remark}[subsection]{Remark}

\title [Flattening and Subanalytic Sets]{Flattening and Subanalytic Sets in Rigid Analytic Geometry}
\author{T.S.~Gardener}
\address {Mathematical Institute\\
University of Oxford\\ 
24-29 St. Giles\\ 
Oxford OX1 3LB (United Kingdom)} 
\email {gardener@maths.ox.ac.uk} 
\author{Hans Schoutens}
\address{Department of Mathematics\\ Rutgers University\\ Hill Center-Bush Campus\\ Piscataway, NJ 08854}
\email{hschoute@math.rutgers.edu}

\begin{document}

\begin{abstract}   
Let $K$ be an algebraically closed field endowed with a complete non-ar\-chi\-me\-dean norm with valuation ring  $R$. Let $f\colon  Y\to X$ be a map of $K$-affinoid varieties. In this paper we  study the analytic structure of the image  $f(Y)\subset X$; such an image is a typical example of a  subanalytic set. We show that  the subanalytic sets are precisely the $\mathbf D$-semianalytic sets, where $\mathbf D$ is the truncated division function first introduced by \name{Denef} and \name{van den Dries}. This result is most conveniently stated as a \QE\ result for the valuation ring $R$ in an analytic expansion of the language of valued fields.

To prove this we establish a Flattening Theorem for affinoid varieties in the style of  \name{Hironaka}, which allows a reduction to the study of subanalytic sets arising from flat maps, that is, we show that a map of affinoid varieties can be rendered flat by using only finitely many local blowing ups. The case of a flat map is then dealt with by a small extension of a result of  \name{Raynaud} and \name{Gruson} showing that the image of a flat map of affinoid varieties is open in the Grothendieck topology. 

Using Embedded Resolution of Singularities, we derive in  the zero characteristic case a Uniformization Theorem for subanalytic sets:  a subanalytic set can be rendered semianalytic using only finitely many local blowing ups with smooth centres. As a  corollary we obtain that any subanalytic set in the plane $R^2$ is semianalytic.  
\end{abstract}

\maketitle


\section{Introduction}

 Subanalytic sets arise naturally in real analytic geometry as  images of proper analytic maps.  The structure of such an image can be quite complicated: it is not necessarily definable by means of inequalities between analytic functions, that is, it is not in general \emph{semianalytic}. Therefore, a compact subset of a real analytic manifold is called \emph{subanalytic}, if it is (at least locally) the projection of a relatively compact semianalytic set. It is then a non-trivial fact that the complement of a subanalytic set is  subanalytic again. Nevertheless, subanalytic sets share many tameness properties with for instance semialgebraic or analytic sets: local finiteness of the number of connected components (which are subanalytic again); the distance between two disjoint closed subanalytic sets is strictly positive; various \Loj\ inequalities hold. Real subanalytic sets were first introduced by \name{\Loj} and subsequently studied by \name{Gabrielov} in \cite{Ga68}{} and \name{Hironaka} in \cite{Hi73a,Hi73c,Hi73b}  by complex analytic methods (flattening, \VE) and by geometric techniques (Resolution of Singularities). A new approach appeared in a paper \cite{DvdD} by \name{Denef} and \name{van den Dries}, where a model-theoretic point of view was taken. This resulted in a much more concise formulation which
has the enormous advantage of being applicable in the $p$-adic context. In this seminal paper, all the basic results in the real case were reproved together with their $p$-adic analogs. 

Motivated by problems of elliptic curves, \name{Tate} constructed a theory of rigid analytic geometry over complete non-archimedean algebraically closed fields. This theory was further developed by \name{Kiehl}, \name{Grauert} et al., largely in analogy with complex analysis. A little later \name{Raynaud} gave an alternative treatment through formal schemes and more recently still, \name{Berkovich} approached the subject from the viewpoint of spectral theory.

From the point of view of model theory, non-archimedean fields provide a very fruitful study. After the work of \name{Ax} and \name{Kochen}, the theory of $\mathbb Q_p$ (the $p$-adics) has been exhaustively studied by \name{Macintyre} and for algebraically closed valued fields long before by \name{Abraham Robinson}. With the recent massive application of model theory to real analytic geometry through the work of \name{van den Dries}, \name{Wilkie} et al., the time seemed right to add analytic structure to the complete algebraically closed valued fields.

Such a study was initiated by \name{Lipshitz}, later on joined by \name{Robinson}, who developed in \cite{Lip93,LR92,LRDim,LRAst} a  theory, allowing more general functions than  rigid analytic ones in the description of semianalytic and subanalytic sets,  thus obtaining a theory of \emph{weak rigid subanalytic sets}. This yielded some important results on rigid subanalytic sets as well (see \ref{2.4}--\ref{2.6} for a further discussion). At the same time the second author obtained a different theory in \cite{Sch94a,Sch94c,Sch94b} were he used a restricted class of analytic functions,  yielding the theory of \emph{strong rigid subanalytic sets}. Both of these theories were based on the model-theoretic approach introduced by \cite{DvdD}. Unfortunately, the same method seems not to work for the general case of a rigid subanalytic set.

It was the insight of \name{Denef} that the methods of \name{Hironaka} might be used in the rigid case as well. The key observation is a result due to \name{Raynaud} and \name{Gruson} describing the image of a flat map between affinoid varieties; this serves as a replacement for the Fibre Cutting Lemma of flat maps in \name{Hironaka}'s work. To make the reduction to the flat case, one needs a good theory of rigid analytic flatificators (to be used as centres of local blowing ups) and the construction of the \VE\ (a compact Hausdorff space encoding finite sequences of local blowing ups; this is the rigid analytic version of the Zariski-Riemann manifold). The former is carried out by the second author in \cite{SchRAF} and the latter by the first author in \cite{GarVE}. However, in order to make the construction of the \VE, it seems necessary to add extra points to the rigid analytic variety, following \name{Berkovich}. The present paper will put all these results together to obtain the sought for theory of rigid subanalytic sets.  Our main theorem states that any rigid subanalytic set can be described by inequalities among functions which are obtained by composition and division of analytic functions (see Section 3 for details). Using the theorem of the complement \cite[Theorem 7.3.2]{LRAst}, it is enough to show this for images of rigid analytic maps. Therefore, we need to \emph{flatten} an arbitrary analytic map by means of local blowing ups. Blowing ups are the cause for having to introduce division in the description of a subanalytic set.  

This flattening procedure is explained in Section 2. A local form in the Berkovich category is derived first from which then a global rigid analytic flattening theorem  is deduced.  In the proof of the former result, we briefly recall some concepts and results from \cite{GarVE} and \cite{SchRAF}. The next Section then contains our main result, preceded by a discussion of the link between blowing ups and the truncated division operator $\mathbf D$. In Section 4 we show how using our main result together with \ERS, one derives a uniformization theorem for rigid subanalytic sets. This is then used to show that subanalytic sets in the plane are in fact semianalytic. The treatment here is analogous to the one in \cite{DvdD}, \cite{Sch94b}, \cite{Sch94b} or \cite{LR92}. In the final section, we have gathered some material on the result of \name{Raynaud} and \name{Gruson}. Most of this is well known, but we needed a small extension of the result, which required an adaptation in the proof as it appeared in \cite{Meh81}. For this reason, it has been added here.

\begin{remark}\label{r:1} 
Rigid analytic geometry can be done over any complete non-archi\-me\-dean field but it is most convenient (and for us at times essential) to assume that the field $K$ is also algebraically closed.  We will make this hypothesis so that the unit disk in $K^{N}$ can be identified with the maximal spectrum of $\scp KS$ where $S=\rij SN$ are independent variables. This unit disk (sometimes also denoted by $\mathbb B^{N}(K)$) is most naturally the $N$-th Cartesian power of the valuation ring $R$ of $K$.  All results of Section 2 remain true for arbitrary complete non-archimedean fields, however to prove a \QE\ Theorem we require the algebraically closed hypothesis.
\end{remark}

\begin{remark}\label{r:2}
The authors have restricted their attention only to the rigid analytic case, but a treatment of Berkovich subanalytic sets seems now to be accessible, using the same methods.
\footnote{To this end, the uniform version \cite[Theorem 7.3.2]{LRAst} of the model-completeness result of \name{Lipshitz} and \name{Robinson} is required.}  
Such a theory would be desirable since then  topological properties of subanalytic sets can be studied, such as the behaviour of connected components, triangulation or even homotopic invariants.
\end{remark}

\begin{remark}\label{r:3}
As far as the characteristic of the field $K$ is concerned, no assumption is needed, except in Section 4 where an application of \ERS\ is used. If one would have a version of \ERS\ in positive characteristic, or, at least of its corollary mentioned in the proof of Theorem~\ref{3.1}, the assumption on the characteristic could be removed.
\end{remark}

 \subsection*{Acknowledgment} 
The authors want to thank \name{Jan Denef}, for sharing with them his idea that \name{Raynaud}'s Theorem could be used in conjunction with a rigid analog of \name{Hironaka}'s work. They are also very grateful to \name{Gabriel Carlyle} for his minute proof reading.

A short survey of the material presented in this paper will appear in \cite{GSQESur}.
 
\section{Rigid Analytic Flattening}

\begin{definition}[Blowing Up]\label{1.1} 
Let $X$ be a \rav.  We will be concerned in this section with local blowing up maps and their compositions. For the definition and elementary properties of rigid analytic blowing up maps, we refer to \cite{SchBU}. Suffice it to say here that they are characterised by the universal property whereby a coherent sheaf of ideals is made invertible. Any blowing up map is proper and an isomorphism away from the centre. If its centre is  nowhere dense, then the blowing up is also surjective. A local blowing up $\pi$ of $X$ is a composition of a blowing up map $\pi'\colon \tilde X\to U$ and an open immersion $U\into X$. We will always assume that $U$ is affinoid.
\footnote{One can always make such an assumption by perhaps shrinking the admissible open, since whenever we use local blowing ups, we will only be interested in a local situation.}

If   $Z$ is the centre of the blowing up $\pi'$ (and hence in particular a closed analytic subvariety of $U$), then we call $Z$ also the centre of $\pi$ and we will say that $\pi$ is the local blowing up of $X$ with locally closed centre $Z$.

 Let $f\colon Y \to X$ be a map of \ravs\ and let $\pi\colon \tilde X\to U\into X$ be a local blowing up with centre $Z\subset U$. If $\theta\colon \tilde Y\to \inverse fU\into Y$ denotes the local blowing up of $Y$ with (locally closed) centre $\inverse fZ$, then by universal property of blowing up, there exists a unique map $\tilde f\colon \tilde Y\to \tilde X$, making the following diagram commute
	\commdiagram[] {\tilde Y} {\theta} {Y} {\tilde f} {f} {\tilde X} 	{\pi} {X.} 
This unique map $\tilde f$ is called the \emph{strict transform} of $f$ under $\pi$ and the above diagram will be referred to as the \emph{diagram of the strict transform}.

In general, we will not be able to work with just a single local blowing up, but we will make use of maps which are finite compositions of local blowing up maps. Therefore, if $\pi\colon \tilde X\to X$ is the composite map  $\psi_1\after\cdots \after\psi_m$, with each  $\psi_{i+1}\colon X_{i+1}\to U_{i}\into X_{i}$  a local blowing up map with centre $Z_i$,  for $i<m$, (with $X=X_0$ and $\tilde X=X_m$), then we define recursively $f_i\colon Y_i\to X_i$ as the strict transform of $f_{i-1}$ under $\psi_i$ where $f_0=f$ and $Y_0=Y$.  The last strict transform $f_m$ is  called the \emph{(final) strict transform} of $f$ under $\pi$ and the other strict transforms $f_i$, for $i<m$, will be referred to as the \emph{intermediate} strict transforms.

For us, the following three possible properties of a map $\pi$ as above, will be crucial.   
\begin{enumerate}  
\item\label{i} The centres $Z_i$ are nowhere dense.  \item\label{ii} The intermediate strict transforms are flat over their centre, that is to say, the restriction $\inverse {f_i}{Z_i}\to Z_i$ is flat, for $i<m$.  \item\label{iii} The final strict transform $\tilde f$ of $f$ under $\pi$ is flat.  \end{enumerate}  Our Flattening Theorem  states that given a map of affinoid varieties $f\colon Y\to X$, we can find finitely many maps $\pi_1, \dots,\pi_s$ as above with these three properties \ref{i}--\ref{iii}, such that, furthermore, the union of their images contains $\op{Im}f$.

\name{Hironaka}'s proof of the complex analytic Flattening Theorem, heavily exploits the fact that complex spaces are Hausdorff and locally compact. The canonical topology (that is, the topology induced by the norm) on a \rav\ is, in general, not locally compact. The Grothendieck topology is not even a topology at all. This means that new ideas are necessary to prove a rigid analytic Flattening Theorem. However, the work of \name{Berkovich} provides us with new analytic spaces, equivalent to rigid analytic spaces as far as their sheaf theory is concerned, but admitting a locally compact Hausdorff topology. We briefly recall their construction. Let $A$ be an affinoid algebra and $X=\aff A$ the corresponding affinoid variety. We fix once and for all a complete normed field $\Upsilon$ extending $K$ and of cardinality big enough so that it contains any completion of any normed extension field of $K$ which is topologically of finite type over $K$. An \emph{analytic point} $x$ of $X$ is defined to be a continuous $K$-algebra morphism $x\colon A\to \Upsilon$.  Let $U=\aff C$ be an affinoid subdomain of $X$ containing $x$, then $U$ is called an \emph{affinoid neighbourhood} of $x$, if the map $x\colon A\to\Upsilon$ factors through a map $C\to\Upsilon$. Two analytic points are said to be congruent if they admit the same system of neighbourhoods. The \emph{affinoid Berkovich space} $\Berk X$ associated to $X$ is then the collection of all congruence classes of analytic points of $X$. Put the weakest topology on $\Berk X$ making all maps $x\mapsto\norm{x(f)}$ continuous, for $f\in A$. This turns the space $\Berk X$ into a compact space, which is  Hausdorff, if $X$ is reduced. As a special case of an analytic point, take any $x\in X$, which then corresponds to some maximal ideal $\maxim$ of $A$. The composite map $A\onto A/\maxim=K\into\Upsilon$ yields an analytic point, called a \emph{geometric point} of $X$. Any affinoid subdomain containing a geometric point is a neighbourhood of that point, and hence $X$ viewed as the set of geometric points can be identified with a subspace of $\Berk X$, and as such is everywhere dense in it.

Finally, one can put a $K$-analytic structure on $\Berk X$, by defining a structure sheaf $\loc_{\Berk X}$ on it. The category of coherent sheaves on $X$ is then equivalent with the category of coherent sheaves on $\Berk X$. So far, we have only introduced the affinoid models in the Berkovich setting. Global models and global morphisms do exist. We refer the reader for further details to \cite{Ber90}, \cite{Ber93} or\cite{SP95}.
\end{definition}

In the statement and the proof of the following Theorem, we will use the term \emph{space} to indicate a Hausdorff, paracompact strictly $K$-analytic Berkovich  space. A \emph{map} between two spaces will be an analytic map in the Berkovich sense. Topological notions are taken with respect to the Berkovich topology.  In particular, a \emph{local blowing up} will be the composition of a blowing up map followed by an open immersion of an affinoid Berkovich space. We always use  black board bold letters $\Berk X, \Berk Y,\dots$ to denote Berkovich spaces.

\begin{theorem}[Local Flattening of Berkovich Spaces]\label{1.2} 
Let $f\colon \Berk Y\to  \Berk X$ be a map of  Hausdorff, paracompact strictly $K$-analytic Berkovich spaces with $\Berk X$ reduced. Pick  $x\in\operate{Im}f$ and let $\Berk L$ be a non-empty compact subset of  $\inverse fx$. There exists a finite collection $E$ of maps $\pi\colon \Berk X_\pi\to\Berk X$, with each $\Berk X_\pi$ affinoid, such that the following four properties hold, where we put  $\Berk X_0=\Berk X$,  $\Berk Y_0=\Berk Y$ and $f_0=f$ and where $\pi\in E$ is arbitrary in the first three conditions.
\begin{enumerate} 
\item\label{1.2.i} The map $\pi$ is a composition $\psi_{0}\after\cdots\after\psi_{m}$ of  finitely many local blowing up maps $\psi_{i}\colon  \Berk X_{i}\to \Berk X_{i-1}$  with locally closed nowhere dense centre $\Berk Z_{i}\subset\Berk X_{i}$, for $i=\range 1 m$. 
\item \label{1.2.ii} Let $f_{i}$ be defined inductively as the strict transform of $f_{i-1}$ under  the local  blowing up $\psi_{i}$. Then $\inverse{f_{i}}{\Berk Z_{i}}\to\Berk Z_{i}$ is flat, for $i=\range 1m$.
\item\label{1.2.iii} The final strict transform $f_m\colon \Berk Y_m\to \Berk X_m$ of $f$  under the whole map $\pi$  given by the strict transform diagram
	\commdiagram[] { \Berk Y_m} {\theta} {\Berk Y } {f_m} f { 	\Berk X_m} {\pi} {\Berk X,} 
is flat at each point of $\Berk Y_m$ lying above a point of  $\Berk L$.
\item\label{1.2.iv}  The union of all the $\op{Im}\pi$, for $\pi\in E$, is a neighbourhood of $x$.
\end{enumerate}
\end{theorem}
\begin{proof} 
\setcounter{stp}0
\begin{step} 
Our first task is to define the \VE\ of an arbitrary space $\Berk X$. The details of this process are in \cite{GarVE}, but the method is wholly due to \name{Hironaka} who makes the construction for complex analytic spaces. 

Let $\mathcal E(\Berk X)$ denote the class of all  maps $\pi\colon \Berk X'\to\Berk X$ which are compositions of finitely many local blowing up maps.  One can define a partial pre-order relation on $\mathcal E(\Berk X)$ by calling $\psi\colon \Berk X''\to\Berk X$ smaller than $\pi$, if $\psi$ factors as $\pi q$, for some $q\colon \Berk X''\to\Berk X'$. We denote this by $\psi\leq\pi$.\footnote{If we want to make this into a partial order, we should take the quotient of $\mathcal E(\Berk X)$ modulo the following equivalence relation: two maps $\pi\colon \Berk X'\to \Berk X$ and $\psi\colon \Berk X''\to \Berk X$ are equivalent if, and only if, $\psi\leq \pi$ and $\pi\leq \psi$. Moreover, the collection of equivalence classes now forms a set, which we denote again by $\mathcal E(\Berk X)$.} The map $q$ is necessarily unique and must belong to $\mathcal E(\Berk X')$ (\cite[Proposition 3.2]{GarVE}). If, moreover, the image $q(\Berk X'')$ of $q$ is relatively compact (that is to say, the closure of its image is compact), then we denote this by $\psi\refine \pi$. Any two maps $\pi_1,\pi_2\in\mathcal E(\Berk X)$ admit a unique minimum or \emph{meet} $\pi_3\in\mathcal E(\Berk X)$ with respect to the order $\leq$ (\cite[Lemma 3.3]{GarVE}), denoted by $\pi_1\en\pi_2$. This meet $\pi_3$ is just the strict transform of $\pi_2$ under $\pi_1$ (or vice versa). With these definitions, $\mathcal E(\Berk X)$ becomes a semi-lattice with smallest element the empty map $\emptyset\colon \emptyset\to\Berk X$. A subset $e$ of $\mathcal E(\Berk X)$ is called a \emph{filter}, if it does not contain $\emptyset$, is closed under meets, and, for any $\psi\in e$ and $\pi\in \mathcal E(\Berk X)$, with $\psi\leq\pi$, we have that also $\pi\in e$. An \emph{\etoile} $e$ on $X$ is now defined as a filter on the semi-lattice  $\mathcal E(\Berk X)$ which is maximal among all filters $e'$ satisfying the extra condition that for any $\pi\in e'$ we can find $\psi\in e'$, with $\psi\refine\pi$.

The collection of all \etoile{}s on $\Berk X$ is called the \VE\ of $\Berk X$ and is denoted by $\mathcal E_{\Berk X}$. This space is topologised by taking for  opens the sets of the form $\mathcal E_\pi$ given  as the collection of all \etoile{}s on $\Berk X$ containing $\pi\colon \Berk X'\to\Berk X$, for some  $\pi\in\mathcal E(\Berk X)$. In fact, $\mathcal E_\pi$ is isomorphic with $\mathcal E_{\Berk X'}$ via the map $J_\pi\colon \mathcal E_{\Berk X'}\to \mathcal E_{\Berk X}$, sending $e'\in\mathcal E_{\Berk X'}$ to the collection of all $\theta\in\mathcal E(\Berk X)$ for which there exists some $\psi\in e'$ such that $\pi\after\psi\leq\theta$, see (\cite[Proposition 3.6]{GarVE}).  The \VE\ is  Hausdorff in this topology (\cite[Theorem 3.11]{GarVE}). Moreover, for any \etoile\ $e\in\mathcal E_{\Berk X}$, the intersection of all $\op{Im}\pi$, where $\pi$ runs through the maps in $e$, is a singleton $\{x\}$ and any open immersion $\restrict 1{\Berk U}\colon \Berk U\into\Berk X$ with $x\in\Berk U$, belongs to $e$ (\cite[Proposition 3.9]{GarVE}). We denote the thus defined map $e\mapsto x$ by  $p_{\Berk X}\colon \mathcal E_{\Berk X}\to\Berk X$. It is a  continuous and surjective map. It is a highly non-trivial result that this map is also proper in the sense that the inverse image of a compact is compact (\cite[Theorem 3.13]{GarVE}).
\end{step}

\begin{step}
Next, we will introduce the concept of a flatificator. Let $f\colon \Berk Y\to\Berk X$ be a map and let $x\in\Berk X$. A \emph{flatificator} of $f$ at $x$ is a locally closed subspace $\Berk Z$ of $\Berk X$ containing $x$, such that $f$ is flat over it (that is to say, the restriction $\inverse f{\Berk Z}\to\Berk Z$ is flat), and such that, whenever $\Berk V$ is a second locally closed subspace containing $x$ over which $f$ is also flat, at least on an open neighbourhood around $x$, then $\Berk V$ is a subspace of $\Berk Z$ locally around $x$ (that is to say, when restricted to some open neighbourhood of $x$). In other words, a flatificator is a largest locally closed subspace over which $f$ becomes flat in a neighbourhood of $x$. Such a flatificator is called \emph{universal}, if it is stable under base change. With this we mean that, if $g\colon \Berk X'\to\Berk X$ is arbitrary, then $\inverse g{\Berk Z}$ is the flatificator of the base change $\Berk Y\times_{\Berk X}\Berk X'\to \Berk X'$ at $x'$, for any $x'$ in the fibre above $x$. In \cite[Theorem A.2]{SchRAF} it is shown that any map $f\colon \Berk Y\to \Berk X$ admits a universal flatificator $\Berk Z$ in each point $x$ of $\op{Im}f$. If $\Berk X$ is moreover  reduced, then we can detect flatness via the flatificator: blowing up the flatificator exhibits some non-trivial portion of non-flatness as torsion. More precisely, it is shown in (\cite[Theorem A.6]{SchRAF}) that  whenever $f$ is not flat in some point of the fibre $\inverse fx$, then there exists a nowhere dense subspace $\Berk Z_0$ of $\Berk Z$, such that the local blowing up $\psi_1\colon \Berk X_1\to\Berk X$ with centre $\Berk Z_0$ \emph{renders the fibre above $x$ smaller}. With this, we mean the following. Let
	\commdiagram[]{\Berk  Y_1 } {\zeta_1}  {\Berk Y} {f_1} f {\Berk  X_1} {\psi_1} {\Berk X} denote the strict transform diagram induced by $\psi_1$. Then, for every $y\in\Berk X_1$  lying above $x$, we have a non-trivial embedding of closed subspaces
	\begin{equation}\label{dagger}
	\inverse{f_1}{y}\times_K\Upsilon \varsubsetneq \inverse 	fx\times_K\Upsilon,
	\end{equation}
where $\Upsilon$ is our universal domain. (Note the extension of scalars is necessary in order to compare these two fibres as subspaces of $\Berk Y\times_K\Upsilon$.)  We refer to this result as the \emph{Fibre Lemma}.
\end{step}

\begin{step}
Fix an \etoile\ $e$ on $\Berk X$ such that $p_{\Berk X}(e)=x$. Put  $x_0=x$, $e_0=e$ and put $f_0\colon \Berk Y_0\to\Berk X_0$ equal to the original map $f\colon \Berk Y\to\Berk X$. Finally, to comply with the enumeration below, put $\Berk Z_{-1}=\Berk X_{-1}=\emptyset$ and let $\psi_0\colon \Berk X_0\to\Berk X_{-1}$ be the empty map. The Fibre Lemma will enable us to define, by induction on $i$,  points $x_i\in\Berk X_i$ and \etoile{}s $e_i\in\mathcal E_{\Berk X_i}$ with $p_{\Berk X_i}(e_i)=x_i$, local blowing up maps $\psi_i\colon \Berk X_i\to\Berk X_{i-1}$ with nowhere dense centre $\Berk Z_{i-1}$, maps $f_i\colon \Berk Y_i\to\Berk X_i$ and non-empty compact subsets $\Berk L_i$ of the fibre $\inverse{f_i}{x_i}$,  at least as long as $f_i$ is not flat in some point of $\Berk L_i$. Each $f_i$ will be the strict transform
	\commdiagram []{ \Berk Y_i} {\zeta_i} {\Berk Y_{i-1} } {f_i} 	{f_{i-1}} { \Berk X_i} {\psi_i} {\Berk X_{i-1}.} 
of the previous map $f_{i-1}$ under the local blowing up $\psi_i$. Moreover, each $f_i$ will be flat above the centre $\Berk Z_i$ and have the property on the fibres \eqref{dagger} in the point $x_{i-1}$. Let us show how to define from the point $x_{i-1}\in\Berk X_{i-1}$ a new point $x_i\in\Berk X_i$ and a  new \etoile\ $e_i$ on $\Berk X_i$. Apply the Fibre Lemma to the point $x_{i-1}$ and the map $f_{i-1}$ to obtain a nowhere dense center $\Berk Z_{i-1}$ and a local blowing up $\psi_i$ rendering the fibre above $x_{i-1}$ smaller as explained in \eqref{dagger}. Since $\Berk Z_{i-1}$ is nowhere dense and contains $x_{i-1}$, we deduce, from \cite[Corollary 3.10]{GarVE}, that $\psi_i\in e_{i-1}$, that is, $e_{i-1}\in\mathcal E_{\psi_i}$. The isomorphism $J_{\psi_i}\colon \mathcal E_{\Berk X_i}\to\mathcal E_{\psi_i}$ then yields a uniquely determined \etoile\ $e_i$ on $\Berk X_i$ and this in turns uniquely determines the point $x_i=p_{\Berk X_i}(e_i)$ of $\Berk X_i$. By diagram chasing, one checks that this implies  $\psi_i(x_i)=x_{i-1}$. Finally, we  define a compact subset of $\inverse{f_i}{x_i}$ by
	\begin{equation}
	\Berk L_i= \Berk Y_i \cap (\{x_i\}\times_{\Berk X_{i-1}} \Berk 	L_{i-1}).
	\end{equation}
Let $\mathcal I_i$ denote the coherent ideal of $\Berk Y\times_K\Upsilon$ defining the fibre $f_i^{-1}(x_i)\times_K\Upsilon$, so that by \eqref{dagger} the chain
	\begin{equation}
	\mathcal I_0\varsubsetneq \mathcal I_1 \varsubsetneq 	\mathcal I_2 \varsubsetneq \dots
	\end{equation}
is strictly increasing on the compact set $\Berk L\times_K\Upsilon$. Therefore this chain must become stationary, say at level $m$, meaning that $f_m$ is flat in each point of $\Berk L_m$. 
 Flatness is an open condition in the source,
\footnote{This is for instance proven in \cite[Corollary A.3]{SchRAF} as a consequence of the existence of a flatificator.}  so that we can find an open $\Berk V$ of $\Berk  Y_m$ containing $\Berk L_m$, such that $\restrict{f_m}{ \Berk V}\colon  \Berk  V\to \Berk X_m$ is flat. Let $\Berk M =  \inverse\theta{ \Berk L}\setminus \Berk V $, where $\theta$ is  the compositum $\zeta_m\after\cdots\after\zeta_1$. Note that $\Berk M$ is compact, since $\theta$ is proper. We claim that after some further local blowing up (in fact, an open immersion will suffice), we may assume that $\Berk M$ is empty, so that the new strict transform will be flat at each point lying above a point of $\Berk L$. To this end, suppose $\Berk M$ is non-empty and pick some $y\in\Berk M$. Since $f_m(y)\neq x_m$, we can  find a compact neighbourhood $\Berk K_y$ of $x_m$ in $ \Berk X_m$,  such that $y\notin\inverse{f_m}{ \Berk K_y}$. Hence
	\begin{equation}\label{}
	\emptyset= \bigcap_{y\in \Berk M} (\inverse{f_m}{ \Berk K_y} 	\cap \Berk M).
	\end{equation}
The compactness of each $\inverse{f_m}{ \Berk K_y} \cap \Berk M $ means that   already a finite number of them, say $\inverse{f_m}{ \Berk K_{y_i}} \cap \Berk M$, for $i<t$, have empty intersection. Let $\Berk K$ be the intersection of these finitely  many $\Berk K_{y_i}$, which is then still a compact neighbourhood of  $x_m$, with the property that
	\begin{equation}
	\inverse{f_m}{ \Berk K} \cap \inverse\theta{ \Berk L} \subset 	\Berk V.
	\end{equation}
Let $\Berk X_e$ be an open of $\Berk X_m$ containing $x_m$ and contained in $\Berk K$. Then the restriction of $f_m$ above $\Berk X_e$ has now the property that it is flat  in each point lying above a point of $\Berk L$. 
\end{step}

\begin{step}
 Summarising, we found for each \etoile\ $e$ on $\Berk X$ with $p_{\Berk X}(e)=x$ a local blowing up map $\pi_e\colon \Berk X_e\to\Berk X$, such that the strict transform
	\commdiagram []{ \Berk Y_e} {\theta_e} {\Berk Y } {f_e} f {	\Berk X_e} {\pi_e} {\Berk X} 
has the property that it is flat in each point lying above a point of $\Berk L$. Moreover, there is a canonically defined   point $x_e$ on $\Berk X_e$ lying above $x$. Let $\Berk C_e$ be a relatively compact open neighbourhood of $x_e$ and set 
	\begin{equation}
	\Berk D_e=J_{\pi_e}(\inverse{p_{\Berk X_e}}{\Berk C_e}).
	\end{equation}
Then, for each $e\in \inverse{p_{\Berk X}}x$, the set $\Berk D_e$ is an open neighbourhood of $e$ and the union of all the $\Berk D_e$ is a neighbourhood of $\inverse{p_{\Berk X}}x$. Let $\{ \Berk H_\lambda\}_{\lambda\in \Lambda }$,  be the  collection of all non-empty relatively compact open neighbourhoods of $x$ in $\Berk  X$. Pick  arbitrary $\lambda_0,\lambda_1\in\Lambda$ with $\overline{\Berk H}_{\lambda_1}\subset \Berk H_{\lambda_0}$, where, in general, $\overline{\Berk H_\lambda}$ denotes the topological closure of  $\Berk H_\lambda$. Consider the open  covering of $\inverse{p_{ \Berk X}}{ \overline{\Berk H}_{\lambda_1}}$ given by  the sets
	\begin{equation}\label{tag2}
	\set { \inverse{p_{ \Berk X}}{ \Berk H_{\lambda_0}} \cap  	\Berk D_e} {e\in \inverse{p_{ \Berk X}}x} \cup \set {\inverse 	{p_{ \Berk X}} {\Berk H_{\lambda_0}\setminus \overline{\Berk  	H}_\lambda}}{\lambda\in\Lambda}.
	\end{equation}
This is indeed a covering, since by the Hausdorff property, we can find  for each $y\neq x$ a $\lambda\in\Lambda$ for which $y\notin \overline{\Berk H}_\lambda$. Since $p_{\Berk  X}$ is proper, we have that $\inverse{p_ {\Berk X}} {\overline{\Berk H}_{\lambda_1}}$ is compact and  hence there exists a finite subset $E\subset\inverse  {p_{\Berk X}}x$ and a finite subset $\Gamma\subset\Lambda$, such  that the collection of all sets of \eqref{tag2} with $e\in E$ and  $\lambda\in\Gamma$ remains a covering of $\inverse{p_{\Berk  X}}{\overline{\Berk H}_{\lambda_1}}$. Putting $\Berk H$ equal to the intersection  of the $\Berk H_\lambda$, for $\lambda\in\Gamma$, this is still a  neighbourhood of $x$ and 
	\begin{equation}
	\inverse {p_{ \Berk  X}}{ \Berk H} \subset \bigcup_{e\in E} 	\Berk D_e.
	\end{equation}
Observing  that $p_{ \Berk X}( \Berk D_e)=\pi_e( \Berk C_e)$ and using that $p_{ \Berk X}$ and $p_{\Berk X_e}$ are surjective, we deduce that
	\begin{equation}
	\Berk H \subset \bigcup_{e\in E} \pi_e( \Berk C_e),
	\end{equation}
as required.
\end{step}
\end{proof}

Apart from the sources already mentioned, a detailed and self-contained explanation of all the terms and properties  used in the above proof can be found in \cite{SchRSS}.

From now on, we  will work in the more familiar category of \ravs\ and consequently we must translate this flattening theorem into a version appropriate for the context. This also calls for a more global result.

\begin{theorem}[Flattening Theorem]\label{1.3} Let $f\colon Y\to X$ be a map of affinoid  varieties with $X$ reduced. Then there exists a finite collection $E$ of  maps $\pi\colon X_\pi\to X$, with each $X_\pi$ again affinoid such that the  following properties hold.
\begin{enumerate} 
\item\label{1.3.i} Each $\pi\in E$ is the composition $\psi_1\after\cdots\after\psi_m$ of finitely many local blowing up maps $\psi_{i}$ with locally closed nowhere dense centre $Z_{i-1}$, for $i=\range 1 m$.
\item\label{1.3.ii}  For each $\pi\in E$, let $f_{i}$ be inductively defined as the strict transform of $f_{i-1}$ under  the local  blowing up $\psi_{i}$. Then $\inverse{f_{i}}{Z_{i}}\to Z_{i}$ is flat, for $i=\range 1m$. The diagram of strict transform is
	\commdiagram []{Y_{i}} {\zeta_{i}} {Y_{i-1}} {f_{i}} {f_{i-1}} 	{X_{i}} {\psi_{i}} {X_{i-1}.}
\item\label{1.3.iii} The strict transform $f_\pi\colon  Y_\pi\to X_\pi$ of $f$ under the whole map $\pi $ (which is $f_m$ according to our enumeration) is flat. The diagram of strict transform is
	\commdiagram []{Y_\pi} \theta Y {f_\pi} f {X_\pi} \pi {X.}
\item\label{1.3.iv} The union of all the $\operate{Im}\pi$, for $\pi\in E$, contains the image
$\op{Im}f$.
\end{enumerate}
\end{theorem}
\begin{proof} 
Let $\Berk X$ and $\Berk Y$ be  the corresponding Berkovich spaces of $X$ and $Y$ respectively and let us continue to write $f$ for  the corresponding map $\Berk Y\to \Berk X$. Fix an analytic point $x$ of  $X$ (that is to say, a point of $\Berk X$), contained in the image of $f$. Let $\Berk L=\inverse fx$, which is closed in $\Berk Y$ whence compact since $\Berk Y$  is. By Theorem~\ref{1.2}, we can find a finite collection $E_x$ of maps $\pi\colon \Berk X_\pi\to\Berk X$ with $\Berk X_\pi$ affinoid,  such that the conditions~\ref{1.2.i}--\ref{1.2.iv} hold. For  each $\pi\in E_x$, let
	\commdiagram []{\Berk Y_\pi} {\theta} { \Berk Y} {f_\pi}f { 	\Berk X_\pi} {\pi} { \Berk X}
be the corresponding strict transform diagram.

By \ref{1.2.iii} of Theorem~\ref{1.2} we  have that the strict transform $f_\pi$ is flat in each point of  $\inverse{\theta}{\inverse fx}=\inverse {f_\pi}{\inverse{\pi}x}$. Let us first show that we can modify the data in such way that   $f_\pi$ becomes flat everywhere.  Since  flatness is open in the source by \cite[Theorem 3.8]{SchRAF}, we can find an open  neighbourhood $\Berk V'$ of $\inverse {f_\pi}{\inverse{\pi}x}$ in  $\Berk Y_\pi$ over which $f_\pi$ is flat. Since $\Berk X_\pi$ and $\Berk Y_\pi$ are compact Hausdorff spaces, we can find  an open neighbourhood $\Berk U'$ of $\inverse{\pi}x$, such that  $\inverse{f_\pi}{ \Berk U'}\subset \Berk V'$. Similarly,  we can find an open  neighbourhood $\Berk U$ of $x$ in $\Berk X$, such that $\inverse{\pi}{ \Berk   U}\subset \Berk U'$.   The   neighbourhood $ \Berk U$ can be taken inside the union of all the  $\operate{Im}{\pi}$, for all $\pi\in E_x$. Set $\Berk  U_\pi= \inverse\pi{\Berk U}$. Note that  $\Berk U_\pi\into  \Berk X_\pi$ is the strict transform of the open immersion $\Berk U\into  \Berk X$ under $\pi$. Let $\psi$ be the restriction of $\pi$ to  $\Berk U_\pi$. The strict transform of $f$ under $\psi$ is the map
	\begin{equation}
	\inverse {f_\pi} { \Berk U_\pi} \to \Berk U_\pi,
	\end{equation}
which by construction is flat, since 
	\begin{equation}
	\inverse {f_\pi} { \Berk U_\pi}\subset \inverse {f_\pi} { \Berk 	U'}\subset \Berk V'.
	\end{equation}
This establishes our claim upon replacing $\pi$ by $\psi$.

Hence we may assume that $f_\pi$ is flat. Note also that in the above process, we have not violated condition~\ref{1.2.iv} of Theorem~\ref{1.2}, so that the $\pi(\Berk X_\pi)$, for all $\pi\in E_x$, form a covering  of an affinoid neighbourhood $\Berk W_x$ of $x$ in $\Berk X$.  We can translate all these diagrams to the rigid  analytic setup and assume that the same diagrams hold with the spaces  now \ravs\  (see Remark~\ref{r:1a} below), where we  keep the same names for our spaces and maps, but just  replace any  blackboard letter, such as $\Berk X,\dots$, by its corresponding roman  equivalent $X,\dots$, denoting the corresponding rigid analytic variety.   In particular, \ref{1.3.i}--\ref{1.3.iii}  hold and we show how to obtain \ref{1.3.iv}.   
 
Let us now vary the analytic point $x$ over  $\op{Im} f$, so that the $W_x$ cover all analytic points of  $\op{Im}f$. Since $\Berk Y$ is compact in the Berkovich topology so is $f(\Berk Y)$. Therefore,  by \cite[Lemma 1.6.2]{Ber93}, already finitely many of the $W_x$ cover all analytic points of $\op{Im}f$. In particular, there is a  finite collection $S$ of analytic points, such that the union of all  $\operate{Im}{\pi}$, for all $\pi\in E_x$ and all $x\in S$, cover $\op{Im}f$. In other words,  condition~\ref{1.3.iv} is now verified as well. 
\end{proof}

\begin{remark}\label{r:1a}
In this translation process from Berkovich data to rigid analytic data, one needs the following. Let $X$ be an arbitrary \rav\ which is quasi-separated with a finite admissible affinoid cover.  Let us denote by  $\Berk X=\Berk  M(X)$  the corresponding Berkovich space. Suppose  $\pi\colon \tilde {\Berk X} \to \Berk U\into \Berk X$ is a local blowing up  with centre $\Berk Z$, where the latter is a closed subspace of the open $\Berk U$. We can  find a wide affinoid $V$ of $X$, such that its closure $\Berk M(V)$ in  $\Berk X$ is contained inside $\Berk U$. Hence there exists a closed  analytic subvariety $Z$ of $V$, such that $\Berk M(Z)=\Berk Z\cap \Berk  M(V)$. Let $p\colon \tilde X \to V$ be the blowing up of $V$ with this centre  $Z$, then $\Berk M(\tilde X)\subset \tilde { \Berk X}$ (see \cite[Lemma 2.2]{GarVE} for the details). So in our  translation we will replace $\pi$ by the (rigid analytic) local blowing  up $\tilde X\to V\into X$. Moreover, if $\Berk W$ is an open inside $\Berk U$ such that its closure $\overline{ \Berk W}$ is still contained in $\Berk U$, then we can take  $V$ such that $\Berk W\subset \Berk M(V)$ and hence 
	\begin{equation}
	\inverse\pi{\Berk W} \subset \Berk M(\tilde X) \subset 	\tilde{\Berk X}.
	\end{equation}
Note that the local blowing up $\tilde{\Berk W}\to\Berk W\into\Berk X$ of $\Berk X$ with centre $\Berk Z\cap\Berk W$ coincides with the restriction $\inverse\pi{\Berk W} \to \Berk X$, so that the rigid analytic local blowing up $\tilde X\to X$ is sandwiched by the Berkovich local blowing ups $\inverse\pi{\Berk W} \to \Berk X$ and $\tilde{\Berk X} \to \Berk X$. The picture is
	\begin{equation}
		\begin{CD}
		\tilde{\Berk W}	@>>>	\Berk W	@>>>	\Berk X		\\
		@VVV		@VVV		@|		\\
	\Berk M(\tilde X)	@>>>	\Berk M(V)	@>>>	\Berk X		\\
		@VVV		@VVV		@|		\\
	\tilde {\Berk X}	@>>>	\Berk U	@>>>	\Berk X,	
		\end{CD}
	\end{equation}
where the composite vertical maps are open immersions and the outer composite horizontal maps are local blowing ups.

Moreover, in this way we can maintain  in the rigid
analytic version all covering properties which were already  satisfied in the Berkovich version.
\end{remark}

\begin{remark}\label{r:2a} Note that we proved something stronger than condition~\ref{1.3.iv}: the union of the images of all $\pi\in E$ covers not only all geometric points of $\op{Im}f$, but also all analytic points.
\end{remark} 
 
\section{Subanalytic Sets}

\begin{definition}\label{2.1}
We now introduce the notion of semianalytic  and subanalytic sets in rigid analytic geometry. There are essentially two different ways of viewing  these objects, one is  geometrical in nature and the other is model-theoretic. We require both of these viewpoints in our analysis. In what follows, let $X=\aff A$ be a reduced affinoid variety (that is to say, $A$ has no non-trivial nilpotent elements).
\end{definition} 
 
\subsection{The Geometric Point of View}

Let $X$ be an affinoid variety.  The subset $ \Sigma \subset X$ is called \emph{semianalytic}, if it is the finite union of basic sets. A \emph{basic subset} is any subset of the form 
	\begin{equation}\label{tag1}
	\set{x\in X}{\regsize{\norm{p_i(x)}\leq \norm{q_i(x)}, \text{ 	for $i<n$ and } \norm{p_i(x)}<\norm{q_i(x)}, \text{ for $n\leq 	i<m$}}}, 
	\end{equation}
with  the $p_i,q_i\in A$. For a \rav\ $X$ which has a finite admissible affinoid covering, we say that $ \Sigma \subset X$ is semianalytic, if the intersection with each element of the cover is semianalytic. Note that if $Y \subset X$ is an inclusion of affinoids, then a set $\Sigma \subset Y$ is semianalytic in $Y$ if, and only if, $\Sigma$ is also semianalytic in $X$.

Once more for $X$ affinoid, we say that $\Sigma \subset X$ is
\emph{subanalytic} if there exists a semianalytic $ \Omega \subset X \times R^{N}$ so that $ \Sigma$ is the image of $\Omega$ under the projection $X\times R^{N} \to X$.  For a \rav\ $X$ with a finite
affinoid cover we say that $ \Sigma \subset X$ is subanalytic, if the intersection with each element of the cover is subanalytic.

Whereas the collection of all semianalytic subsets of $X$ is easily seen to be a Boolean algebra, this is no  longer obvious at all for the class of subanalytic sets. Recently, \name{Lipshitz} and \name{Robinson} gave a proof of this result in \cite[Corollary 1.6]{LRAst}. Below, we give a short review of their results, since we will make use of them in the proof of our Quantifier  Elimination~\ref{2.7}.

In order to give a neat description of a subanalytic set, it is convenient to introduce a special function $\mathbf D$, first introduced by \name{Denef}  and \name{van den Dries} in their paper \cite{DvdD}, in which they describe  $p$-adic subanalytic sets. Put 
	\begin{equation}\label{eq:D}
	\mathbf D\colon R^2\to R\colon  (a,b)\mapsto 
	\begin{cases}
	a/b\qquad\qquad&\text{if $\norm a \leq\norm b\neq0$}\\ 
	0\qquad\qquad&\text{otherwise.}
	\end{cases}
	\end{equation}
Let us denote by $\scp AS$ with $S=\rij SN$, the ring of \emph{strictly convergent power series} over $A$, where we endow $A$ with its supremum norm. We define  the \emph{algebra $A^{\mathbf D}$ of $\mathbf D$-functions}  on  $X$, as the smallest $K$-algebra of $K$-valued functions on $X$ containing $A$ and closed under the following two operations. 
\begin{enumerate}  
\item\label{D1} If $p,q\in A^{\mathbf D}$, then also $\mathbf D(p,q)\in A^{\mathbf D}$. 
 \item\label{D2} If $p\in\scp A{T_1,\dots,T_N}$ and $q_i\in A^{\mathbf D}$ with $\norm{q_i}\leq 1$, for $i=\range 1N$, then also $p\rij  qN\in A^{\mathbf D}$. 
\end{enumerate}
Here, the function  $\mathbf D(p,q)$ is to be considered as  a pointwise division, that is, defined by $x\mapsto \mathbf D(p(x),q(x))$.  Note also that if $p\in A^{\mathbf D}$ then $p$ defines a bounded function on $X$ and hence it makes sense to define $\norm  p=\sup_{x\in X}\norm{p(x)}$. If we allow in the definition of  semianalytic sets also $\mathbf D$-functions rather  than just elements of $A$, we may now formulate the definition of  \emph{$\mathbf D$-semianalytic} sets: the functions appearing in \eqref{tag1} may be elements of $A^{\mathbf D}$.  Our main result says that a set is $\mathbf D$-semianalytic if, and only if, it is subanalytic. 
 
\subsection{The Model-Theoretic Point of View}

The basic model will be the valuation ring $R$ which we think of as the domain of convergence for the functions in $\scp KX$ with $X$ a single variable. More generally $R^{N}$ is the unit ball ${\mathbb B}^{N}(K)$ which is the domain of convergence for the functions in the Tate ring in $N$ variables.
 
 The \emph{analytic language}  $\anlang$ for $R$ consists of two $2$-ary relation symbols $\mathbf P_\leq$ and $\mathbf P_<$ and  an $n$-ary function symbol $F_f$, for every strictly convergent power series $f$ in $n$-variables of norm at most one, that is to say,  for every $f\in\scp R{X_1,\dots,X_n}$, where $n=0,1,\dots$. The interpretation of $R$ as an $\anlang$-structure is as  follows. Each $n$-ary function symbol $F_f$ is interpreted as the corresponding function $f\colon  R^n\to R$, defined by the  strictly convergent power series $f$ (note that $\norm f\leq1$, so that $f$ is indeed $R$-valued). The relation symbol $\mathbf  P_\leq$ interprets the subset $\set {(x,y)\in R^2} {\norm x\leq\norm y}$ of $R^2$, and likewise, $\mathbf P_<$ describes the  subset $\set {(x,y)\in R^2} {\norm x<\norm y}$. Hence, the atomic formulae in this language (or rather, their interpretation  in $R$) are of the following three types  
	\begin{subequations}
	\begin{gather}
	f(x)=g(x),\label{taga}\\ 
	\norm{f(x)}\leq \norm{g(x)},\label{tagb}\\
	\norm{f(x)} < \norm{g(x)}.\label{tagc} 
	\end{gather}
	\end{subequations}
Note that the first type can be rewritten as $\norm{f(x)-g(x)}\leq 0$, so that we actually only have to deal with types  \eqref{tagb} and \eqref{tagc}. One can of course define $\mathbf P_<(x,y)$ as $\niet\mathbf P_\leq(y,x)$, but the  advantage of not doing so is that all formulae can now be made equivalent with positive ones, that is, without using the  negation symbol. In this language the semianalytic sets in $R^{N}$ correspond to the quantifier-free definable sets; the subanalytic sets are existentially definable. There is an example, adapted from \name{Osgood}'s example, of a subanalytic set which is not semianalytic.  This means that it is not the case that every formula in the language is equivalent in the first order $\anlang$-theory of $R$ to a quantifier free formula. (Note that a similar failure of quantifier elimination is observed in the real  or in the $p$-adic case. )

To remedy this, we introduce an expansion  $\anDlang$ of $\anlang$ with one new $2$-ary function symbol $\mathbf D$,  which we will interpret in our structure as the function $\mathbf D$ given by \eqref{eq:D}. If $K$ were the $p$-adic field (and hence  $R=\mathbb Z_p$), then  by a theorem of \name{Denef} and \name{van  den Dries} \cite{DvdD}, $R$ admits Elimination of Quantifiers in  an expansion of this language where one needs to add extra predicates, one for each $n=2,3,\dots$, to express that an element is an $n$-th power; a similar expansion occurs in \name{Macintyre}'s algebraic Quantifier Elimination for $\mathbb Z_p$. In the algebraically closed case, these predicates are unnecessary. Hence the following is the  natural rigid analytic analogue: the valuation ring $R$ of $K$ admits  Elimination of Quantifiers in the language $\anDlang$. 
 
 Let us see how this ties in with the notion of subanalyticity. A subset of $R^N$ which is definable in the language  $\anlang$ by a quantifier free formula, is precisely a semianalytic set whereas an existentially definable set is   a projection of a semianalytic set, so consequently  subanalytic. It is not too hard to see that the function $\mathbf D$ is existentially definable and whence also  every $\mathbf D$-function on $R^N$, so that any $\mathbf D$-semianalytic subset of $R^N$ is subanalytic. The claim that every subanalytic set is $\mathbf D$-semianalytic is then equivalent with Quantifier Elimination in the language $\anDlang$.
\footnote{By induction on the number of quantifiers, it is enough to eliminate only existential quantifiers to obtain Quantifier Elimination.}

Let $X$ be an arbitrary affinoid variety and choose some closed immersion $j\colon X \into R^{N}$ for some $N \in \nat$. Hence $j(X)$ is quantifier-free definable in the language $\anlang$ and $X$ is isomorphic with $j(X)$.  More generally, any quasi-compact rigid analytic variety is isomorphic with some quantifier-free  $\anlang$-definable set.  Also note that semianalytic sets (respectively, subanalytic sets) in such a variety $X$ correspond to quantifier-free definable (respectively, existentially definable) subsets of $X$ under the chosen closed immersion.

We will be adopting from now on the geometric point of view. In particular,  we will identify $\aff(\scp K{S_1,\dots,S_n})$ with $R^n$. 

\begin{example}\label{2.1.3} 
If $f\colon Y\to X$ is a map of affinoid varieties, then the image $f(Y)$ is a typical subanalytic subset of $X$ (not necessarily semianalytic!). Subanalyticity follows from projecting the graph of $f$ (which is  analytic, whence semianalytic) onto $X$. More generally, it follows that $f(\Sigma)\subset X$ is subanalytic whenever $\Sigma\subset Y$ is subanalytic. This example shows that even when one is merely interested in closed analytic  subsets, one needs to study subanalytic sets as well. However, there are some particular kinds of maps which have better understood image. For instance, \name{Kiehl}'s Proper Mapping Theorem \cite{Ki67} (or \cite[9.6.3. Proposition 3]{BGR}) states that the image of a proper map is closed analytic. However, this does not tell us anything about the image of a semianalytic set under a proper map. In fact, in \cite{Sch94a,Sch94b} the second author shows that if $\Sigma\subset Y$ is semianalytic and $f\colon Y\to X$ is proper, then $f(\Sigma)$ is $\mathbf D$-semianalytic in $X$; he carries out a systematic study of the sets arising in this way--the \emph{strongly subanalytic} sets. One might hope though that certain proper maps, viz. blowing up maps,  nevertheless behave better with respect to semianalyticity. In Proposition~\ref{2.3} below, as a first step towards Quantifier Elimination, we will show that the image of a semianalytic (and even of a  $\mathbf D$-semianalytic) under a blowing up map is $\mathbf D$-semianalytic, at least away from the centre. It is because of this (rather straightforward) result that $\mathbf D$-functions enter in the description of subanalytic sets. We note that for surfaces we have a much stronger statement.  It
follows from \cite[Proposition 3.1]{Sch94c} that if $\pi\colon X \to R^2$ is the blow up of the plane in a point and if $ \Sigma \subset X$ is semianalytic, then $ \pi (\Sigma)$ is also semianalytic.  This we use in the proof of Corollary~\ref{3.2}.

A second class of affinoid maps with well-understood images are the flat maps: their images are finite unions of rational domains and hence in particular semianalytic. This highly non-trivial result is due to \name{Raynaud} and \name{Gruson} (a full account by \name{Mehlmann} appeared in \cite{Meh81}).  Because of its crucial role in our argument and since we need a slight improvement of their original result in the form Theorem~\ref{2.2} below, we will provide most of the details in Section~\ref{app}. 
\end{example}

\begin{theorem}[Raynaud-Gruson-Mehlmann]\label{2.2} 
Let $f\colon  Y\to X$ be a flat map of affinoid varieties. Let $\Sigma$ be a  semianalytic subset of $Y$ given by
	\begin{equation}
	\Sigma=\set{y\in Y} {\norm{h_i(y)}<1\text{ and }\norm{g_j(y)}\geq 1,  \text{ for } i\in I\text{ and } j\in J}
	\end{equation}
where $ h_{i} , g_{j}$ are finitely many functions in the affinoid algebra of $Y$ each of whose supremum norm is at most one. Then $f(\Sigma)$ is semianalytic in $X$.  
\end{theorem}
\begin{proof} See Section~\ref{app}.
\end{proof}
 
\begin{proposition}\label{2.3}  
Let $\pi\colon \tilde X\to X$ be a map of rigid analytic varieties and let $\Sigma$ be a $\mathbf  D$-semianalytic subset of $\tilde X$. If $\pi$ is a locally closed immersion, then $\pi(\Sigma)$ is $\mathbf D$-semianalytic in  $X$. If $\pi$ is a local blowing up map with centre $Z$, then $\pi(\Sigma)\setminus Z$ is $\mathbf D$-semianalytic in $X$. 
 \end{proposition} 
\begin{proof}  
For closed immersions, the statement is trivial. If $U=\aff C\into X=\aff A$ is a rational  affinoid subdomain,  then $C=\scp A{f/g}$, where $f=\rij fn$ with $f_i,g\in A$ having no common zero. Hence any function $h\in C$ defined on  $U$ is $\mathbf D$-definable on $X$ (just replace any occurrence of $f_i/g$ by $\mathbf D(f_i,g)$). Now, any affinoid subdomain  is a finite union of rational subdomains by \cite[7.3.5. Corollary 3]{BGR} and hence we proved the proposition for any  affinoid open immersion as well. From this, the general locally closed immersion case follows easily.  
 
 This leaves us with the case of a blowing up. Without loss of generality, we may assume $X$ to be affinoid. Let us briefly  recall the construction of a blowing up map as described in \cite{SchBU}. Let $X=\aff A$ and let $Z$ be a closed analytic  subvariety of $X$ defined by the ideal $\rij g n$ of $A$.  We can represent $A$ as a quotient of some $\scp KS$, with  $S=\rij Sm$, so that $X$ becomes a closed analytic subvariety of $R^m$.  However, in order to construct the blowing up of  $X$ with centre $Z$, we need a different embedding, given by the surjective algebra morphism
	\begin{equation}
	\scp K{S,T} \onto A \colon T_j\mapsto g_j, 
	\end{equation}
for  $j=\range 1n$, extending the surjection $\scp KS\onto A$ and where $T=\rij Tn$. This gives us a closed immersion  $i\colon X\into R^m\times R^n$ and after identifying $X$ with its image $i(X)$, we see that $Z=X\cap (R^m\times 0)$. Now, the  blowing up $\pi\colon\tilde X\to X$ is given by a strict transform diagram  
	\commdiagram []{\tilde X} {\pi} {X} {\tilde \imath} {i} {W} 	{\gamma} {R^m\times R^n} 
where  $\gamma$ denotes the blowing up of $R^m\times R^n$ with centre the linear space $R^m\times 0$. There is a  standard finite admissible affinoid covering $\{W_1, \dots, W_n\}$ of $W$ where each $W_j$ has affinoid algebra  
	\begin{equation}
	C_j= \frac {\scp K{S,T,U}}{(T_jU_1-T_1,\dots,T_jU_n-T_n)}, 
	\end{equation}
so  that $\gamma(s,t,u)=(s,t)$ for any point $(s,t,u)\in W_j$, where the latter is considered as a closed analytic subset of  $R^m\times R^n\times R^n$ via the above representation of $C_j$. Moreover, $\tilde X$ is a closed  analytic subvariety of $X\times_{(R^m\times R^n)}W$. Therefore, if we set $\tilde X_j=\inverse{\tilde \imath}{W_j}$, then $\{\tilde  X_1,\dots,\tilde X_n\}$ is a finite admissible affinoid covering of $\tilde X$ with the affinoid algebra $\tilde A_j$ of each  $\tilde X_j$ some quotient of the affinoid algebra  
	\begin{equation}\label{eq:aff}
	\frac {\scp A{\hat U_j}}{(g_jU_1-g_1,\dots,g_jU_n-g_n)}
	\end{equation}
of $W_j \times_{(R^m\times R^n)}X$, where $\hat U_j$ means all variables $U_k$ save $U_j$. 
 
 With this notation, let us return to the proof of the proposition. We are given some $\mathbf D$-semianalytic set $\Sigma$ of  $\tilde X$ and we seek to describe the image $\pi(\Sigma)\setminus Z$. Let us focus for the time being at one  $\Sigma\cap\tilde X_j$, where $j\in\{1,\dots,n\}$. Since $\Sigma\cap \tilde X_j$ is $\mathbf D$-semianalytic, we can find a  quantifier free $\anDlang$-formula $\varphi_j(\tuple s,\tuple u)$, such that $(s,u)\in R^m\times R^n$ belongs to $\Sigma\cap  \tilde X_j$ if, and only if, $\varphi_j(s,u)$ holds.  Hence, for $s\in R^m$, we have that    
	\begin{equation}\label{eq:form}
	s\in\pi(\Sigma\cap \tilde X_j), \qquad\text{if and  only 	if,}\qquad(\exists\tuple u) \varphi_j(s,\tuple u).
	\end{equation}
If $\varphi_j(s,u)$ holds, then in particular $(s,u)\in\tilde X_j$ and hence by \eqref{eq:aff}, we have that $g_j(s)u_k=g_k(s)$, for all  $k=\range 1n$. Now, a point $s\in R^m$ does not belong to $Z$, precisely when one of the $g_k(s)$ does not vanish.  Therefore, as $j$ ranges through the set $\{1,\dots,n\}$ and using \eqref{eq:form}, it is not too hard to see that $s\in R^m$ belongs to $\pi(\Sigma)\setminus Z$ if, and only if, $s\in  X$ and 
	\begin{equation*}
	\Of_{j=1}^n \En_{k=1}^n \norm{g_k(s)}\leq\norm {g_j(s)} \en 	g_j(s)\neq0 \en \varphi_j(s,\mathbf 	D(g_1(s),g_j(s)),\dots,\mathbf D(g_n(s),g_j(s))),
	\end{equation*}
which is indeed a $\mathbf D$-semianalytic description of $\pi(\Sigma)\setminus Z$. 
 \end{proof} 
 
\begin{remark}  
The above result is unsatisfactory in so far as it does not tell us anything about $\pi(\Sigma)$  restricted  to the centre $Z$ of the blowing up. If we could prove that also $\pi(\Sigma)\cap Z$ were $\mathbf D$-semianalytic, then the  whole image $\pi(\Sigma)$ would be $\mathbf D$-semianalytic, as we would very much like to show. However, above $Z$,  the map $\pi$ looks like a projection map, so that we can't say much more about $\pi(\Sigma\cap\inverse\pi Z)=\pi(\Sigma)\cap Z$ except that it is a subanalytic set.  If $Z$ would be zero dimensional and whence finite, then clearly also $\pi(\Sigma)\cap Z$ is $\mathbf D$-semianalytic. This  suggests that we might be able to use the above result in order to prove Quantifier Elimination by an induction  argument  on the dimension of $X$, as soon as we can arrange that $Z$ has strictly smaller dimension than $X$. This will be  the case, if $Z$ is nowhere dense; a condition we ensure will always be fulfilled.  
 
 Another point ought to be mentioned here: although a blowing up $\pi\colon \tilde X\to X$ is an isomorphism outside its  centre $Z$, this does \emph{not} automatically imply that one can deduce from the $\mathbf D$-semianalyticity of  $\Sigma\setminus\inverse\pi Z$ the same property for its (isomorphic) image $\pi(\Sigma)\setminus Z$. What is going on here is  that being ($\mathbf D$-)semianalytic is not an intrinsic property of a set, but of its embedding in a larger space. In other words, being  isomorphic as point sets is not enough and thus the above statement is not a void one.   
\end{remark} 

Before we turn to the proof of our main theorem, let us give a brief review on the model-completeness result of \name{Lipshitz} and \name{Robinson}. Geometrically, this amounts to the fact that the complement of a subanalytic set is again subanalytic. This is by no means a straightforward  result. In the real case it was originally proved by \name{Gabrielov} in \cite{Ga68} using quite involved arguments, which later became simplified by \name{Bierstone} and \name{Milman} in \cite{BM88}. In the paper \cite{DvdD} of \name{Denef} and \name{van den Dries} an entirely different approach was taken, using a much stronger result, namely, the class of  subanalytic sets is equal to the class of  $\mathbf D$-semianalytic sets. The complement of a $\mathbf D$-semianalytic set is evidently
$\mathbf D$-semianalytic.

\begin{theorem}[Lipshitz-Robinson]\label{2.4} 
The complement $X\setminus\Sigma$ and the closure $\overline\Sigma$ (in the canonical topology) of a subanalytic set $\Sigma$ in $X$, where $X$ is a reduced quasi-compact \rav, is again subanalytic.
\end{theorem}

\begin{theorem}[Lipshitz-Robinson]\label{2.5} 
Let $X$ be a reduced quasi-compact \rav\ and let $\Sigma$ be a subanalytic set in $X$. Then there exists a finite partition of $\Sigma$ by pairwise disjoint rigid analytic submanifolds $X_i$ of $X$ such that their underlying sets are subanalytic in $X$.
\end{theorem}

The proofs of both Theorems rely on a certain \QE\ result in some appropriate language and we refer the reader to the papers \cite[Theorem 7.4.2]{LRAst} and  \cite[Theorem 4.4]{LRDim} by \name{Lipshitz} and \name{Robinson}. Let us just show how  one can derive a good dimension theory for subanalytic sets from these results. First, there is the  notion of the dimension of a quasi-compact \rav. This is defined as the  maximum of the Krull dimensions of all its local rings (we give the empty  space dimension $-\infty$). In case $X=\aff A$ is affinoid, this is just the Krull dimension of $A$. Next, we define the dimension of a subanalytic set  $\Sigma$ in $X$ as the maximum of all $\op{dim} Y$, where  $Y\subset\Sigma$ is a submanifold of $X$. If $\Sigma$ carries already the  structure of a manifold, then clearly its subanalytic dimension equals its  manifold dimension.

The relevant properties for this dimension function are now summarised by the following proposition.

\begin{proposition}\label{2.6} 
Let $X$ be a quasi-compact \rav\ and let $\Sigma$ and $\Sigma'$ be (non-empty) subanalytic sets in $X$. Then the following holds.
\begin{enumerate}
\item\label{2.6.i} If $\Sigma\subset\Sigma'$, then the dimension of $\Sigma$ is at most the dimension of $\Sigma'$.
\item\label{2.6.ii} The dimension of $\Sigma$ is zero if, and only if, $\Sigma$ is finite.
\item\label{2.6.iii} The dimension of $\Sigma$ equals the dimension of its closure (in the canonical topology) $\overline\Sigma$.
\item\label{2.6.iv} The dimension of the \emph{boundary} $\overline\Sigma\setminus\Sigma$ is strictly smaller than the dimension of $\Sigma$.
\item\label{2.6.v} If $f\colon X\to Y$ is a map of quasi-compact \ravs, then the dimension of $f(\Sigma)$ is at most the dimension of $\Sigma$, with equality in case $f$ is injective.
\item\label{2.6.vi} If $\Sigma$ is semianalytic, then the dimension of $\Sigma$ is equal to the (usual) dimension of its Zariski closure.
\end{enumerate}
\end{proposition}

\begin{remark} Note that by Theorem~\ref{2.4} both the closure $\overline\Sigma$ and the boundary  $\overline\Sigma \setminus\Sigma$  are indeed subanalytic.
\end{remark}

\begin{proof}  
The  first two statements follow from the fact that the dimension of a  subanalytic set is the maximum of the dimensions of each manifold in any  finite subanalytic manifold partitioning  (as in Theorem~\ref{2.5}). The other  statements require more work. See \cite{Lip93} and, in particular, for Property~\ref{2.6.iv}, see \cite[Theorem 4.3]{LRDim}. See also \cite[3.15-3.26]{DvdD} for the $p$-adic  analogues--the proofs just carry over to our present situation, once one  has Theorem~\ref{2.5}.
\end{proof}

\begin{theorem}[Quantifier Elimination]\label{2.7} 
For any reduced affinoid variety $X$, the subanalytic subsets of $X$  are precisely the $\mathbf D$-semianalytic subsets of $X$. 
\end{theorem} 
\begin{proof}   
We have already seen that $\mathbf D$-semianalytic sets are subanalytic.  To prove the converse, let $\Sigma$ be a subanalytic set of $X$. We will  induct on  the dimension of $\Sigma$ and then on the dimension of $X$.  The zero-dimensional case follows immediately from \ref{2.6.ii} in Proposition~\ref{2.6}. Hence fix  $\op{dim}\Sigma=k>0$ and $\op{dim}X=d>0$. In particular, $k\leq d$.

\setcounter{stp}0
\begin{step}
It suffices to take $\Sigma$ closed in the canonical topology. Indeed, assume the theorem proven for all subanalytic sets which are closed in the canonical topology. Let $\overline\Sigma$ be the closure of $\Sigma$ with respect to the canonical topology. By Theorem~\ref{2.4} and \ref{2.6.iii} of Proposition~\ref{2.6}, also $\overline\Sigma$ is subanalytic and of dimension equal to the dimension of $\Sigma$. Hence by our  assumption $\overline\Sigma$ is even $\mathbf D$-semianalytic. Let $\Gamma$ be the boundary $\overline\Sigma\setminus\Sigma$, which is again subanalytic by Theorem~\ref{2.4}. Moreover, by \ref{2.6.iv} of Proposition~\ref{2.6}, $\Gamma$ has strictly smaller dimension than $\Sigma$. Hence, by our induction hypothesis on the dimension of a subanalytic set, we have that also $\Gamma$ is $\mathbf D$-semianalytic. Therefore also $\Sigma = \overline\Sigma\setminus\Gamma$, as required.
\end{step}

\begin{step}
Hence we may assume that $\Sigma$ is closed in the canonical topology.  There exists a  semianalytic subset $\Omega'\subset X\times R^N$, for some $N$, such that  $\Sigma=f'(\Omega')$, where $f'\colon X\times R^N\to X$ is the projection on the first factor. The union of finitely many $\mathbf  D$-semianalytic sets is again such. Therefore, without loss of generality, we may even take $\Omega'$ to be a basic set, that is to say,  of the form  
	\begin{equation}
	\set{(x,t)\in X\times R^N} {{\En_{i<m}\norm{p_i(x,t)} 
	\leq \norm{q_i(x,t)}\ \en \En_{m\leq i<n} \norm{p_i(x,t)}< 	\norm{q_i(x,t)}}}, 
	\end{equation}
where the $p_i$ and $q_i$ are in $\scp AT$, with $X=\aff A$ and $T=\rij TN$. Introduce $n$ new variables $Z_i$ and  consider the  closed analytic subset $Y$ of $X\times R^{N+n}$ given by the equations $p_i-Z_iq_i=0$, for $i<n$.  Let $\Omega$ be the basic subset of $Y$ given by $(x,t,z)\in Y$ belongs to $\Omega$ whenever $\norm{z_i}<1$, for $m\leq i<n$. Let $q$ be the  product of all the $q_i$, for $m\leq i<n$, and we obviously can assume that $q\neq 0$ if $\Sigma$ is non-empty. If  $f\colon Y\to X$ denotes the composition of the closed immersion $Y\into X\times R^{N+n}$ followed by the projection  $X\times R^{N+n}\to X$, then $f(\Omega\cap U)=\Sigma$, where $U$ is the complement in $Y$ of the  zero-set of $q$.  Using \cite[Corollary 2.2]{Sch97}, we may, after perhaps modifying some of the equations defining $Y$, assume that the closure of $U$ in the canonical topology equals the whole of $Y$ and hence the closure (in the canonical topology) of $\Omega\cap U$ is $\Omega$. Now $\Omega\cap U\subset \inverse f\Sigma$ and so $\Omega=\overline{\Omega\cap U}\subset \inverse f\Sigma$, since $\Sigma$ is closed and $f$ is continuous. Hence $f(\Omega)= \Sigma$.
\end{step}

\begin{interlude}
 Before giving the details of the remaining steps, let's pause to give a brief outline of how we will go about it. According to our  Flattening  Theorem~\ref{1.2}, we can find finitely many diagrams 
	\commdiagram [daggerpi]{Y_\pi} {\theta_\pi} {Y} {f_\pi} {f} 	{X_\pi} {\pi}{X} 
where each $\pi$ is a finite composition of local blowing up maps with the properties \ref{1.2.i}--\ref{1.2.iii} and such that $\op{Im}f$ is contained in the union of all the $\operate{Im}\pi$.  Now, in order to study $\Sigma=f(\Omega)$, we will chase $\Omega$ around  these diagrams~\eqref{daggerpi}. There are only finitely many $\pi$ to consider; it will suffice to do this for one such $\pi$ since the analysis for the others is identical. First we take the preimage $\inverse{\theta_\pi}\Omega$, which is again a semianalytic set defined by inequalities  of the form $\norm {h(s)}<1$ where the $h$ are functions on $Y_\pi$ of supremum norm at most one. Next we take the image  of the latter set under $f_\pi$. Our extension of \name{Raynaud}'s Theorem (Theorem~\ref{2.2}) guarantees that this image is  semianalytic. Finally we push this set back to $X$ via $\pi$ and denote this set temporarily by $\Sigma'$. If we had  the full version of Proposition~\ref{2.3}, namely, that a local blowing up map preserves $\mathbf D$-semianalyticity, then this last set would be  indeed $\mathbf D$-semianalytic. 

Of course, in chasing $\Omega$ around the diagram, we might have lost some points. In other words, it may well be the case that $\Sigma'\neq \Sigma$. However, this could happen only for points coming from one of the centres of the local blowing ups that make up $\pi$ (since outside its  centre, a blowing up map is an isomorphism). Above each of these centres the strict transform is flat so we  account for those missing points using Theorem~\ref{2.2} once more. Hence the only problem in the  this reasoning lies in the application of Proposition~\ref{2.3}: it is not the whole image that we can account for by means of that  proposition, but only for the part outside the centre. Now an induction on dimension allows us to dispose of this part too.
\end{interlude}

\begin{step}
Our second induction hypothesis says that any subanalytic  set in an affinoid variety of dimension  strictly smaller than $d$ is $\mathbf D$-semianalytic.  From this we obtain a stronger version of  Proposition~\ref{2.3}.

\begin{lemma}\label{2.3p}  
Let $\pi\colon \tilde W\to W$ be any local blowing up of a quasi-compact \rav\ $W$ of dimension at most $d$ whose centre $Z$ is nowhere dense. If $\Gamma\subset\tilde W$ is $\mathbf D$-semianalytic, then $\pi(\Gamma)\subset W$ is also $\mathbf D$-semianalytic.
\end{lemma}

The key point is that $Z$ has dimension strictly smaller than the dimension $d$ of $W$, which is also  the dimension of $\tilde W$. We have an equality
	\begin{equation}
	\pi(\Gamma)=(\pi(\Gamma)\setminus Z) \cup 	(Z\cap\pi(\Gamma)).
	\end{equation}
By Proposition~\ref{2.3} we know that $\pi(\Gamma)\setminus Z$ is $\mathbf D$-semianalytic  and by our induction hypothesis on the dimension we also have that $Z\cap\pi(\Gamma)$ is $\mathbf D$-semianalytic in $Z$ (take a finite affinoid covering to reduce to the affinoid case). This means that $ \pi(\Gamma)$ is $D$-semianalytic in $W$.
\end{step}

\begin{step}
We return to the proof of Theorem~\ref{2.7}. According to Theorem~\ref{1.3}, there exists a finite collection  $E$ of maps $\pi\colon  X_\pi\to X$, such that each $\pi\in E$ induces a strict  transform diagram~\eqref{daggerpi} with properties \ref{1.3.i}--\ref{1.3.iv}. (The intermediate strict transform diagrams are given by \eqref{daggeri} below). By \ref{1.3.iv}, if we could show that each $\operate{Im}\pi\cap \Sigma$ is $\mathbf D$-semianalytic in $X$, then  the same would hold for $\Sigma$, since there are only finitely many $\pi$. Therefore, let us concentrate on one such $\pi=\pi_1\after\dots\after\pi_m$ and adopt the notation from \ref{1.1} for this map, so that in particular,  \ref{1.3.i}--\ref{1.3.iii}  holds. Let each $\pi_{i+1}$ be the blowing up of the admissible  affinoid $U_i\subset X_i$ with nowhere dense centre $Z_i\subset U_i$.  The  diagram of strict transform is given by  
	\commdiagram [daggeri]{Y_{i+1}} {\theta_{i+1}} {Y_i} {f_{i+1}} 	{f_i} {X_{i+1}}{\pi_{i+1}} {X_i.}  
Define  inductively $\Omega_i \subset Y_i$ as $\inverse{\theta_i}{\Omega_{i-1}}$ starting from  $\Omega_0=\Omega$. Note that each $\Omega_i$ is a semianalytic set of $Y_i$ defined by several inequalities of the type  $\norm h<1$, where each $h\in\loc(Y_i)$ is of supremum norm at most one.  Define also inductively, but this time by downwards  induction, the sets $W_{i-1}=\pi_i(W_i)\subset U_i\subset X_i$ where we start with $W_m=X_m=X_\pi$. In particular, we have that  $W_0=\operate{Im}\pi$. By \ref{2.3p} each $W_i$ is $\mathbf D$-semianalytic in $X_i$. In order to describe $\Sigma$, we will furthermore make  use of the sets $\Gamma_i$ defined as $f_i(\Omega_i)\cap W_i$, for  $i\leq m$. In  particular, note that $\Gamma_0$ is nothing other than $f(\Omega)\cap W_0=\Sigma\cap\operate{Im}\pi$, which we aim to show is  $\mathbf  D$-semianalytic. 
 
The next claim shows how two successive members in the chain of commutative diagrams~\eqref{daggeri} relate the $\Gamma_i$: for each $i< m$, we have an equality 
	\begin{equation}\label{ddaggeri}
	\Gamma_i = \pi_{i+1}(\Gamma_{i+1}) \cup (\Gamma_i\cap 	Z_i).
	\end{equation}
Assume we have already established \eqref{ddaggeri}, for each $i<m$. We will prove, by downwards induction on $i\leq m$,  that each $\Gamma_i$ is $\mathbf D$-semianalytic in $X_i$, so that in particular $\Gamma_0$ would be $\mathbf  D$-semianalytic in $X$, as required. First, since $f_\pi=f_m$ is assumed to be flat, we can apply Theorem~\ref{2.2}  to $\Omega_m$ to conclude that $\Gamma_m=f_m(\Omega_m)$ is semianalytic whence $\mathbf D$-semianalytic in  $X_m$. Assume now that we have already proved that $\Gamma_{i+1}$ is $\mathbf D$-semianalytic in $X_{i+1}$  and we want to obtain the same conclusion for $\Gamma_i$ in $X_i$. Using \eqref{ddaggeri}, it is enough to establish  this for both sets in the right hand side of that equality. The first of these, $\pi_{i+1}(\Gamma_{i+1})$, is $\mathbf  D$-semianalytic since we have now the strong version~\ref{2.3p} of Proposition~\ref{2.3} at our disposal. As for the second set,  $\Gamma_i\cap Z_i$, also this one is $\mathbf D$-semianalytic, since $f_i$ restricted to $\inverse{f_i}{Z_i}$ is flat and since
	\begin{equation}
	\Gamma_i\cap Z_i=f_i(\Omega_i \cap \inverse{f_i}{Z_i}) \cap 	W_i,
	\end{equation}
so that Theorem~\ref{2.2}  applies. Note that we already established that $W_i$ is $\mathbf D$-semianalytic. 
\end{step}
  
\begin{step}
Therefore,  it only remains to prove \eqref{ddaggeri}. To show that 	\begin{equation}
	\Gamma_{i} \cap Z_{i} \supset f_i(\Omega_i \cap \inverse{f_i}{Z_i}) \cap W_i
	\end{equation}
is straightforward.  We show the reverse inclusion.  Let $x_{i} \in \Gamma_{i}$. That means that there exists $y_i\in\Omega_i$ and $w_{i+1}\in W_{i+1}$ such that  $f_i(y_i)=x_i=\pi_{i+1}(w_{i+1})$. If $x_i\in Z_i$,   we are done.  Hence assume that  $x_i\notin Z_i$ so that $y_i\notin\inverse{f_i}{Z_i}$.   However, since $W_i\subset U_i$ we have  that $y_i\in\inverse{f_i}{U_i}$. Since $\theta_{i+1}$ is the  blowing up of $\inverse{f_i}{U_i}$ with centre  $\inverse{f_i}{Z_i}$ whence an isomorphism outside this centre, we can even find $y_{i+1}\in Y_{i+1}$, such  that $\theta_{i+1}(y_{i+1})=y_i$. From $y_i\in\Omega_i$ it then follows that $y_{i+1}\in\Omega_{i+1}$. Put  $x_{i+1}=f_{i+1}(y_{i+1})$.  Commutativity of the strict transform diagram implies that  $\pi_{i+1}(x_{i+1})=x_i=\pi_{i+1}(w_{i+1})$. Since $x_i\notin Z_i$, the blowing up $\pi_{i+1}$ is injective at that  point, so that $w_{i+1}=x_{i+1}$ which therefore belongs to $f_{i+1}(\Omega_{i+1})\cap W_{i+1}=\Gamma_{i+1}$,   proving our claim, and hence also our main theorem.  
\end{step}
  \end{proof}

\begin{remark}  
We can derive from this proof also a weak uniformization as follows. Define  $\Sigma_i$ inductively as the inverse image of $\Sigma_{i-1}$ under $\pi_i$, for $1\leq i\leq m$, with $\Sigma_0=\Sigma$.  With notations as in the above proof, we can derive, for $i<m$, from \eqref{ddaggeri} the following
identity  
	\begin{equation}
	\Sigma_{i+1}\cap W_{i+1} = \Gamma_{i+1} \cup 	\big(\inverse{\pi_{i+1}}{\Gamma_i\cap W_i} \cap W_{i+1}\big). 
	\end{equation}
For $i=m-1$, this takes the simplified form $\Sigma_m=\Gamma_m\cup\inverse{\pi_m}{\Gamma_{m-1}\cap Z_{m-1}}$. Now, as already  observed,  $\Gamma_m$ is semianalytic in $X_m=X_\pi$ and similarly $\Gamma_{m-1}\cap Z_{m-1}$ is semianalytic in  $X_{m-1}$ and whence also its preimage under $\pi_m$. In other words, we showed the following corollary.  
\end{remark}
 
\begin{corollary}\label{2.8}  
Let $X$ be a reduced affinoid variety and let $\Sigma$ be a subanalytic set  in $X$. There exists a finite collection of compositions of finitely many local blowing up maps $\pi_1,\dots,\pi_n$ with nowhere dense centre, such that the union of the  $\operate{Im}{\pi_i}$ contains $\Sigma$, and such that  each preimage $\inverse {\pi_i}\Sigma$ has become semianalytic.
\end{corollary} 
\begin{proof} 
This follows from the above discussion in the  case where $\Sigma$ is closed in the canonical topology. The reduction to this case uses an induction argument similar to the one in the proof of the theorem. 
\end{proof} 
 
Note also that to prove the corollary, we do not make use of Proposition~\ref{2.3} but only of Theorem~\ref{2.2}. For an improvement of Corollary~\ref{2.8}, at  least in the zero characteristic case, see the Uniformization Theorem~\ref{3.1} below, where we will be able to take  smooth centres for the blowing ups involved. 
 
\section{Uniformization}

In  \cite[Theorem 4.4]{Sch94b} it was proved that for any strongly subanalytic set $\Sigma$ in an affinoid manifold $X$, there  exists a finite covering family of compositions $\pi$ of finitely many local blowing ups with smooth and nowhere dense centre, such that  the preimage $\inverse\pi \Sigma$ is semianalytic, provided  the characteristic of $K$ is zero. The restriction to zero characteristic is entirely due to the lack of an Embedded Resolution  of Singularities in positive characteristic.  \name{Hironaka}'s Embedded Resolution of Singularities for varieties in characteristic zero can be applied to the rigid analytic setting.  See \cite[Theorem 3.2.5]{SchERS} for details.  We are now able to extend the Uniformization Theorem  to the class
of all subanalytic sets.

\begin{theorem}[Uniformization]\label{3.1}  
Let $X$ be an affinoid manifold (that is, an affinoid variety, all local rings of which are regular) and assume  $K$ has characteristic zero. Let $\Sigma$ be a subanalytic subset of $X$. Then there exists a finite collection $E$ of maps $\pi\colon X_\pi\to X$, with each $X_\pi$ again affinoid,  such that the  following properties hold.
\begin{enumerate}
\item \label{3.1.i} Each $\pi\in E$ is the composition $\psi_1\after\cdots\after\psi_m$ of finitely many local blowing up maps $\psi_i$ with nowhere dense and smooth centre, for $i< m$.
\item \label{3.1.ii} The union of all the $\operate{Im}\pi$, for $\pi\in E$, equals $X$.
\item \label{3.1.iii} For each $\pi\in E$, we have that $\inverse\pi\Sigma$ is semianalytic in $X_\pi$.
\end{enumerate}
\end{theorem}
\begin{proof}  
Let $X=\aff A$. We use the following corollary to the Embedded Resolution Theorem: given $p,q\in A$, then there exists a finite collection $E'$ of maps, such that  \ref{3.1.i} and  \ref{3.1.ii} hold, for  each $\pi\colon X_\pi\to X$ in $E'$, and furthermore either $p\after\pi$ divides $q\after\pi$,  or vice versa,  $q\after\pi$ divides $p\after\pi$, in the affinoid  algebra of  $X_\pi$. See for instance \cite[Lemma 4.2]{Sch94b} for a proof. 
 
From our Quantifier  Elimination~\ref{2.7}, we know that $\Sigma$ is  $\mathbf D$-semianalytic. By a (not too difficult)  argument, involving an induction on the number of times the function $\mathbf D$ appears in one of the describing  functions of $\Sigma$ (for details see \cite[Theorem 4.4]{Sch94b}), we can reduce to the case that there is only one such occurrence. In other words, we may assume that  there exist a quantifier free formula $\psi(\tuple x,\pmb y)$  in the language $\anlang$ and functions $p,q\in A$, such that  $x\in \Sigma$ if, and only if,  
	\begin{equation}\label{}
	\psi(x,\mathbf D(p(x),q(x)))\quad\text{holds}.
	\end{equation}
After an appeal to the corollary of Embedded Resolution of Singularities applied to $p$ and $q$, and since we only seek to prove our result modulo finite collections of maps for which \ref{3.1.i} and  \ref{3.1.ii} holds, we  may already assume that either $p$ divides $q$ or $q$ divides $p$. In the former case, there is some $h\in A$, such that  $q=hp$ in $A$. Therefore, $\mathbf D(p(x),q(x))=0$, unless $q(x)\neq0$ and $\norm{h(x)}=1$, in which case it is equal to  $1/h(x)$.  Let $U_1$ be the affinoid subdomain defined by $\norm {h(x)}\leq 1/2$ and $U_2$ by $\norm{h(x)}\geq  1/2$, so that $\{U_1,U_2\}$ is an admissible affinoid covering of $X$. Hence $x\in U_1$ belongs to $\Sigma$ if, and only if,  $\psi(x,0)$ holds, whereas $x\in U_2$ belongs to $\Sigma$ if, and only if,   
	\begin{equation}
	\relax[\norm{h(x)}\geq 1  \en q(x)\neq 0 \en \psi(x,1/h(x))] 	\of [(\norm{h(x)}<1 \of q(x)=0) \en\psi(x,0)] 
	\end{equation}
holds.  Observe that $1/h$ belongs to the affinoid algebra of $U_2$, since $h$ does not vanish on $U_2$. In other  words, $\Sigma$ is semianalytic on both sets and whence on the whole of $X$. 
 
In the remaining  case that $q$ divides $p$, so that there is some $h\in A$, such that  $qh=p$ in $A$, we have an even simpler description of $\Sigma$, namely $x\in \Sigma$, if and  only if, 
	\begin{equation}
	\relax[p(x)\neq0\en\psi(x,h(x))] \of [p(x)=0\en\psi(x,0)] 
	\end{equation}
holds, again showing that $\Sigma$ is semianalytic. 
 \end{proof} 
 
We wish to state a further corollary which is a strengthening of a result in \cite[Theorem 3.2]{Sch94c}.  The proof of \cite[Theorem 3.2]{Sch94c} is applicable directly to the context here (see for instance the final remark in the introduction of \cite{Sch94c}).
 
\begin{corollary}\label{3.2} 
Suppose $K$ has characteristic zero and let $\Sigma\subset R^2$. If $\Sigma$ is subanalytic, then  in fact it is semianalytic.  
\end{corollary} 

\begin{remark} 
Using \name{Abhyankar}'s Embedded Resolution of Singularities \cite{Ab} in positive characteristic for excellent local rings of dimension two, one can remove the assumption on the characteristic in the Corollary.
\end{remark}

\section{Elimination along Flat Maps}\label{app}

This  section will be devoted to a proof of Theorem~\ref{2.2}. In it, we will need some properties of the  \emph{reduction functor} applied to an affinoid algebra. However, for our purposes, we do not need to introduce the whole  machinery of reductions but can make do with an ad hoc construction to be presented below. First, let us fix some notation.  We will denote the maximal ideal of $R$ by $\wp$,  that is, $\wp$ is the collection of all $r\in K$, such that $\norm r<1$. Note that $R$ is a \emph{rank-one} valuation ring, so that the only prime ideals are $(0)$ and $\wp$. The residue field $R/\wp$ will be denoted by $\bar K$. Notice that it is also an  algebraically closed field.    

\begin{definition}\label{A.1}  
We will call  an $R$-algebra $\adm A$  \emph{topologically of finite type}, if $\adm A$ is a homomorphic image of some $\scp  R{S_1,\dots,S_m}$. From  $\adm A$,  we can construct an affinoid algebra $A$ by tensoring with $K$, namely let $A=\adm A\tensor_RK$. If $\adm A$ is moreover flat over $R$, then $\adm A\subset A$.  If we  start with an affinoid  algebra $A=\scp KS/I$ and define $\adm A$ as $\scp RS/\adm I$ where $\adm I=I\cap\scp RS$, then  $\adm A$ is torsion-free whence flat over $R$. By tensoring with $K$ we recover our original affinoid algebra, that is, $A=\adm A\tensor_RK$. However, $\adm A$ depends on the particular choice of representing $A$ as a homomorphic image of some $\scp KS$.

For the sake of simplicity, let us assume that $K$ is algebraically closed.\footnote{This assumption is not essential, although the proofs would require some modifications for the general case; see for instance \cite[Corollary 5.6.10]{SchRSS}.}  Let $\adm A$ be an flat $R$-algebra which is topologically of finite type and let $A=\adm A\tensor_RK$ be the corresponding affinoid algebra. With respect to  the structure map $R\to \adm A$, any prime ideal of  $\adm A$ lies either above $(0)$  or above $\wp$. The prime ideals lying above $(0)$ are in one-one correspondence with  $\op{Spec}A$, whereas the prime ideals lying above $\wp$ are in one-one correspondence with $\op{Spec}\adm A/\wp\adm A$. Therefore, although $\adm A$ is in general not Noetherian, $\op{Spec}\adm A$ has finite combinatorial dimension, as it is the (disjoint) union of $\op{Spec}A$ and $\op{Spec}\adm A/\wp\adm A$ (note that $\adm A/\wp\adm A$ is a finitely generated $\bar K$-algebra). Let us call a morphism $\adm x\colon \op{Spec}R\to \operate{Spec}{\adm A}$ of $R$-schemes an \emph{$R$-rational point}. Then to give a point  $x\in \aff A$ (that is to say, a maximal ideal $\maxim$ of $A$) is the same as to give an $R$-rational point $\adm x$ (given by the ideal $\maxim\cap\adm A$).  The \emph{reduction} of $x$ is the restriction of $\adm x$ to the closed immersion $\op{Spec}\bar K\into\op{Spec}R$ and is denoted by $\bar x$. In other words, $\bar x$ is given by the maximal ideal $(\maxim\cap\adm A)+\wp\adm A$.  Let us denote the  maximal spectrum of $\adm A$ by $\operate{Max}{\adm A}$. The \emph{reduction map} $\xi\colon \aff A\to \operate{Max}{\adm A}$ is the map given by sending $x$ to its reduction $\bar x$. A word of caution:  the reduction map is not  induced by any algebra morphism.

The reduction map $\xi$ is functorial in the following sense. Let  $\adm\varphi\colon \adm A\to\adm B$ be an $R$-algebra morphism with $\adm A$ and $\adm B$ flat and topologically of finite type. Let $\varphi\colon A\to B$ be the  morphism of affinoid algebras obtained by tensoring $\adm \varphi$ with $K$. Then we have a commutative diagram   
	\commdiagram [diamond]{\aff B} {f} {\aff A} {\xi} {\xi} 	{\operate{Max}{\adm B}} {\adm f} {\operate{Max}{\adm A}}
where $f$ and $\adm f$ are  the respective maps on the maximal spectra induced by $\varphi$ and $\adm\varphi$.

It is well known (see for instance \cite{Meh81}) that the reduction map is surjective (regardless whether $K$ is algebraically closed or not). This  is an immediate consequence of the Flat Lifting Lemma below. 
\end{definition}

\begin{lemma}[Flat Lifting Lemma]\label{A.2} 
Let $\adm A\to\adm B$ be an $R$-algebra morphism with $\adm A$ and $\adm B$ $R$-flat and topologically of finite type. Let $f\colon \adm Y=\operate{Spec}{\adm B}\to\adm X=\operate{Spec}{\adm A}$ denote the corresponding map of affine schemes. Let $\adm x\colon \op{Spec}R \to \adm X$, be an $R$-rational point of $\adm X$ and let $\bar x$ denote its reduction $\op{Spec}\bar K\to \adm X$. Suppose there exists a $\bar K$-rational point $\bar y\colon \op{Spec}\bar K \to \adm Y$, such that
	\commtrianglefront [cd1] {\op{Spec}\bar K} {\bar y} {\adm Y} f 	{\adm X} {\bar x}
commutes. If $f$ is flat (that is to say, if $\adm A\to \adm B$ is flat), then there exists an $R$-rational point $\adm y$ of $\adm Y$, which has reduction $\bar y$, and is such that
	\commtrianglefront [cd2]{\op{Spec} R} {\adm y} {\adm Y} f 	{\adm X} {\adm x}
commutes. We call $\adm y$ a \emph{factorisation of $\adm x$ lifting $\bar y$}. 
\end{lemma}
\begin{proof} 
Let $\adm\pr$ be the prime ideal of $\adm A$ associated to $\adm x$. Let $\bar\pr=\adm\pr +\wp\adm A$, so that it is the maximal ideal of $\adm A$ associated to $\bar x$. Finally, let $\bar{\mathfrak q}$ be the maximal ideal of $\adm B$ associated to $\bar y$, so that the commutativity of \eqref{cd1} translates into
	\begin{equation}\label{tag3}
	\bar\pr = \bar{\mathfrak q}\cap\adm A.
	\end{equation}
Since $\adm A\to\adm B$ is flat, the Going Down Theorem (see for instance \cite[Theorem 9.5]{Mats}) guarantees the existence of a prime ideal $\adm{\mathfrak n}$ of $\adm B$, such that $\adm{\mathfrak n}\subset\bar{\mathfrak q}$ and
	\begin{equation}\label{tag4}
	\adm\pr = \adm{\mathfrak n} \cap \adm A.
	\end{equation}
Since $\adm\pr\cap R=(0)$, also $\adm{\mathfrak n}\cap R=(0)$. Among all  prime ideals $Q$ of $\adm B$, with $\adm{\mathfrak n}\subset Q\subset \bar{\mathfrak q}$ and $Q\cap R=(0)$, let $\adm{\mathfrak q}$ be one of maximal height (recall that $\op{Spec}\adm B$ has finite combinatorial dimension). Note that $Q=\adm{\mathfrak n}$ is at least one prime ideal satisfying both conditions. 

We claim that $\adm {\mathfrak q}$ determines an $R$-rational point. Assuming the claim, it follows from $\adm{\mathfrak q}\subset \bar{\mathfrak q}$ that $\adm{\mathfrak q}$ is a  lifting of $\bar{\mathfrak q}$. By Formulae~\eqref{tag3} and \eqref{tag4}, we have inclusions
	\begin{equation}
	\adm\pr \subset \adm{\mathfrak q}\cap \adm A \subset \bar\pr.
	\end{equation}
Since $\adm{\mathfrak q}\cap R=(0)$, we cannot have that $\adm{\mathfrak q}\cap\adm A=\bar\pr$. Since $\adm A/\adm\pr\iso R$, the only two
prime ideals of $\adm A$ containing $\adm\pr$ are $\bar\pr$ and $\adm\pr$ itself. Therefore, we conclude that
	\begin{equation}\label{eq:IV.3.3.7}
	\adm{\mathfrak q}\cap \adm A=\adm\pr,
	\end{equation}
so that the $R$-rational point $\adm {\mathfrak q}$ is a factorization of $\adm\pr$, as required.

It remains to prove that $\adm{\mathfrak q}$ determines an $R$-rational point. In other words, if we set $\adm C=\adm B/\adm{\mathfrak q}$, then we need to show that $\adm C=R$. Write $\adm C$ as a homomorphic image of some $\scp RS$. Since $\adm C$ is a domain and since $\adm{\mathfrak q}\cap R=(0)$, it follows that $\adm C$ has no $R$-torsion and consequently it must be flat over $R$. Put $C=\adm C\tensor_RK$, so that $C$ is an affinoid algebra containing $\adm C$ (here we used flatness). By Noether Normalization (\cite[6.1.2. Theorem 1]{BGR}), we can find a finite injective homomorphism $\phi\colon\scp KT\into C$, for some variables $T=\rij Tn$. Moreover, inspecting the proof of Noether Normalization, one checks that $\phi$ is given by $T_i\mapsto p_i(S)$ with each $p_i$ a polynomial of Gauss norm $1$; see \cite[5.2.4 and 5.1.3. Corollary 7]{BGR} for a discussion. In particular, $\phi$ maps $\scp RT$ into $\adm C$. On the other hand, let $f\in \adm C$ and let $a_i\in \scp KT$ be such that $f^n+\phi(a_1)f^{n-1}+\dots+\phi(a_n)=0$ in $C$. By [BGR 6.2.2 Prop 4], the supremum norm of $f$ is equal to the maximum of the $\norm{a_i}^{1+i}$. Since $f\in \adm C$ its supremum norm is at most $1$. Therefore, all $a_i$ have Gauss norm at most $1$, whence belong to $\scp RT$. In conclusion, the restriction $\scp RT\into \adm C$ is integral. Therefore, $\bar{\mathfrak q}\adm C\cap\scp RT$ is a maximal ideal of $\scp RT$ (\cite[Corollary 4.17]{Eis}). Since $\bar K$ is algebraically closed, it follows that
	\begin{equation}
	\bar{\mathfrak q}\adm C\cap\scp RT=(T_1-x_1,\dots,T_n-x_n)\scp RT+\wp\scp RT
	\end{equation}
for some $x_i\in R$. If $n>0$, then by the Going Down Theorem (\cite[Theorem 13.9]{Eis}) for the normal domain $\scp RT$, we can find a (non-zero) prime ideal $\mathfrak h$ of $\adm C$ contained in $\bar{\mathfrak q}\adm C$ and such that $(T_1-x_1)\scp RT=\mathfrak h\cap\scp RT$. By the maximality of $\adm{\mathfrak q}$, we have that $\mathfrak h\cap R\neq 0$, so that $\wp\adm C\subset\mathfrak h$. Consequently, $\wp\scp RT\subset (T_1-x_1)\scp RT$, contradiction. In other words, we must have that $n=0$, so that $\adm C$ is integral over $R$. Since the field of fractions of $\adm C$ is then algebraic over $K$, it must be equal to it, as we assumed $K$ to be algebraically closed. Therefore  $R\subset \adm C\subset K$. Since $R$ is a valuation ring and clearly $\adm C\neq K$, we conclude that $\adm C=R$, as required.
\end{proof} 

Let $\adm A$ be an flat $R$-algebra, topologically of finite type, with corresponding affinoid algebra $A=\adm A\tensor_RK$. Applying this lemma to the flat map $R\to\adm A$ and the $R$-rational point given by the identity morphism, shows the surjectivity of the reduction map $\xi$. 

The following observation will be constantly used below in the proof of Theorem~\ref{2.2}. Let $h_1,\dots,h_s\in\adm A$ and let $\Sigma$ denote the semianalytic set of all $y\in\aff A$, such that  $\norm{h_i(y)}<1$, for $i<r$ and $\norm{h_i(y)}\geq 1$, for $r\leq i<s$. Call such a set \emph{special}. Let $\adm \Sigma$ denote the  locally closed subset of $\operate{Max}{\adm A}$ consisting of all maximal ideals $\bar \maxim$, such that  $h_i\in\bar\maxim$, for $i<r$, and $h_i\notin\bar\maxim$, for $r\leq i<s$. Using the surjectivity of the reduction map, one easily verifies that these two sets are related to one other by  
	\begin{equation}
	\inverse\xi{\adm\Sigma}=\Sigma \qquad\qquad\text{and} 	\qquad\qquad \adm\Sigma=\xi(\Sigma).
	\end{equation}
In other words, $\xi$ induces a bijection between the class of finite Boolean combinations of special subsets of  $\aff A$ and the class of constructible subsets of $\op{Spec}\bar A$, where $\bar A=\adm A/\wp \adm A$.

\begin{proof}[{Proof of Theorem~\ref{2.2}}] So we are given a flat map $f\colon \aff B\to \aff A$ of affinoid varieties and a  special set $\Sigma$ of $\aff B$.  Suppose that 
	\begin{equation}
	\Sigma =\set {y\in\aff B} {\regsize{\norm{h_i(y)}\diamond_i 1 \text{ for 	$i<s$}}},
	\end{equation}
where $h_i\in B$ are of supremum norm at most one and $\diamond_i$ is either $<$ or $\geq$. We want to prove that $f(\Sigma)$ is semianalytic.

 Using Proposition~\ref{A.4} below, we may reduce to the following case. There exist flat $R$-algebras $\adm A$ and $\adm  B$ which are topologically of finite type, such that $A=\adm A\tensor_RK$ and $B=\adm B\tensor_RK$, and such that $h_i\in\adm B$, for all $i<t$,  and there exists a flat morphism of $R$-algebras $\adm A\to \adm B$  which induces the map $f$ (after tensoring with $K$).   By our observation above, there exists a  locally closed set $\adm \Sigma$ of $\operate{Max}{\adm  B}$ such that $\inverse\xi{\adm\Sigma}=\Sigma$. Since $\operate{Max}{\adm B}$  can be identified with  $\operate{Max}{\bar B}$, where $\bar B=\adm B/\wp\adm B$, we can view $\adm  \Sigma$  as a locally closed subset of the this space as well. If  we also put  $\bar A=\adm A/\wp\adm A$, then since both rings are   finitely generated $\bar K$-algebras, we can invoke \name{Chevalley}'s  Theorem  to conclude that the image of  $\adm\Sigma$ under the induced map $\bar f\colon \operate{Max}{\bar  B}\to\operate{Max}{\bar A}$  is a constructible set $\adm \Omega$.  Identifying  $\operate{Max}{\adm  A}$ with $\operate{Max}{\bar A}$, we  may  consider $\adm\Omega$ as a constructible set of the former space as  well and as such it is the image of $\adm \Sigma$  under the map $\adm  f$ induced by $f$. Let $\Omega=\inverse\xi{\adm\Omega}$, so that by our above observation $\Omega$ is semianalytic in $\aff  A$. Hence, we will have proved our theorem once we show that
	\begin{equation}
	f(\Sigma)=\Omega.
	\end{equation}
The commutative diagram \eqref{diamond} of \ref{A.1} expressing the functoriality of $\xi$, provides the inclusion  $f(\Sigma)\subset \Omega$, so we only need to deal with the opposite inclusion. 
 
To this end, let $x\in\Omega$. Let $\adm x$ be the corresponding $R$-rational point and let $\bar x$ be the reduction $\xi(x)$ of $x$. By assumption, $\bar x\in\adm\Omega$ and hence it is the image under $\adm f$ of some point $\bar y\in\adm\Sigma$. Since $\adm f$ is flat, we can apply  Lemma~\ref{A.2} to  obtain an $R$-rational point $\adm y$ factoring through $\adm x$ and lifting $\bar y$. In other words, if $y\in\aff B$ denotes the point corresponding to $\adm y$, then this translates into $f(y)=x$ and $\xi(y)=\bar y$. Since $\bar y\in\adm \Sigma$, the latter implies that $y\in\Sigma$, as required.    
\end{proof}

\begin{proposition}\label{A.4} 
Let $f\colon Y=\aff B\to X=\aff A$ be a flat map of  affinoid varieties and let $h_j\in B$, for $j<t$, be of supremum norm at most one.  There exist finite coverings $\{U_i=\aff A_i\}_{i<s}$ of $X$ and  $\{V_i=\aff B_i\}_{i<s}$ of $Y$ by rational subdomains and $R$-algebra  morphisms $\varphi_i^\circ\colon  A_i^\circ\to   B_i^\circ$ of flat $R$-algebras which are topologically of finite type, such that, for all $i<s$, we have that
\begin{enumerate}
\item\label{l:1} $A_i=A_i^\circ\tensor_RK$  and $B_i=B_i^\circ\tensor_RK$,
\item\label{l:2} the morphism $\varphi_i^\circ$ is flat and induces the map $\restrict f{V_i}\colon V_i\to U_i$,
\item\label{l:3} $h_j\in B_i^\circ$, for all $j<t$.
\end{enumerate}
\end{proposition}
\begin{proof} 
Since the $h_j$ are of norm at most one and using \cite[6.4.3. Theorem 1]{BGR}, we can find a flat and topologically of finite type $R$-algebra $\adm B$ containing all $h_j$ with $\adm B\tensor_RK=B$, a flat and topologically of finite type $R$-algebra $\adm A$ with $\adm A\tensor_RK=A$ and an $R$-algebra morphism $\adm\varphi\colon \adm A\to\adm B$ inducing the map $f$. In general, $\adm\varphi$ will not be flat. To remedy this, we use \cite[Theorem 3.4.8]{Meh81}, in order to find admissible coverings as asserted, for which \ref{l:1} and \ref{l:2} hold. Moreover, from its proof it follows that $B_i^\circ$ is a homomorphic image of $A_i^\circ\tensor_{\adm A}\adm B$. Therefore, also \ref{l:3} is satisfied.
\end{proof}

\begin{remark} 
The result in \name{Mehlmann}'s paper is quite an intricate matter, using \name{Raynaud}'s approach on rigid analysis through formal schemes and admissible formal blowing ups; an alternative proof can be found in \cite{BL93}.
\end{remark}

\bibliographystyle{amsplain}
\bibliography{myabbrev.bib,references.bib}

\end{document}